\documentclass[10pt,reqno]{amsart}

\usepackage[margin=1in]{geometry} 
\usepackage{amsmath,amsthm,amssymb,mathrsfs}
\usepackage[latin1]{inputenc}
\usepackage{tikz}
\usepackage{url}
\usepackage{graphicx}
\usepackage{multirow}
\usepackage{subfigure}
\usepackage{float}
\usepackage[ruled,vlined]{algorithm2e}
\usepackage[hidelinks]{hyperref}
\usepackage{adjustbox}
\usepackage{booktabs}

\usepackage{lineno} 

\usepackage{pgfplots}
\pgfplotsset{compat=newest}

\usepackage{tikz}
\usetikzlibrary{positioning}

\newtheorem{thm}{Theorem}[section]

\newtheorem{exm}[thm]{Example}

\newtheorem{assu}{Assumption}[section]

\numberwithin{equation}{section}

 \newcommand\blfootnote[1]{%
  \begingroup
  \renewcommand\thefootnote{}\footnote{#1}%
  \addtocounter{footnote}{-1}%
  \endgroup
}

\allowdisplaybreaks
\begin{document}
\title[Two-scale NN Singularly-perturbed Dynamical Syms]{{Two-scale} neural networks for singularly perturbed dynamical systems with multiple parameters}

\author{Qiao Zhuang}
\address[Qiao Zhuang]{School of Science and Engineering, University of Missouri-Kansas City, Kansas City, MO 64110, USA}
\email{qzhuang@umkc.edu}

\author{Taorui Wang}
\address[Taorui Wang]{Department of Mathematical Sciences, Worcester Polytechnic Institute, Worcester, MA 01609, USA}
\email{twang13@wpi.edu}

\author{Rita Wanjiku}
\address[Rita Wanjiku]{School of Science and Engineering, University of Missouri-Kansas City, Kansas City, MO 64110, USA}
\email{rwf66@umkc.edu}

\author{Majid Bani-Yaghoub}
\address[Majid Bani-Yaghoub]{School of Science and Engineering, University of Missouri-Kansas City, Kansas City, MO 64110, USA}
\email{baniyaghoubm@umkc.edu}

\author{Zhongqiang Zhang}
\address[Zhongqiang Zhang]{Department of Mathematical Sciences, Worcester Polytechnic Institute, Worcester, MA 01609, USA}
\email{zzhang7@wpi.edu}

\begin{abstract}
We extend our two-scale neural-network method for scalar singularly perturbed problems with one small parameter 
to dynamical systems 
with multiple small parameters. 
To accommodate multiple small parameters, we use a single effective scale parameter defined as the
geometric mean of all parameters. We thus augment the network input with a scale-aware feature, enabling it to capture sharp solution transitions intrinsically.
 Numerical experiments across a range of dynamical systems demonstrate that 
 the proposed framework can handle coupled systems with multiple and high-contrast small parameters and obtain satisfactory accuracy in capturing solution features induced by small parameters.
\end{abstract}


\date{\today}
\maketitle

\section{Introduction}\label{sec:intro}
%
 In this work, we extend the use of two-scale neural networks (2SNN)
 and successive training \cite{twoscaleNN_2024} within the framework of 
 physics-informed neural networks (PINNs) to the first- and second-order {singularly perturbed dynamical systems} with \textit{multiple small parameters} associated with the leading differential operators. %

%
In \cite{twoscaleNN_2024}, we develop 2SNN for scalar problems with a single parameter, where the scale parameter is straightforward, and so is successive training.  However, for problems with multiple small parameters, it is unclear which parameter should govern scale augmentation or curriculum learning (successive training), or whether we should choose several scale parameters where curriculum learning should be carefully designed.
Such difficulties are further amplified for the singularly perturbed dynamical systems considered in this work, which involve coupled components and multiple fast-slow dynamics with strong scale contrasts.
To preserve simplicity and ease of use in \cite{twoscaleNN_2024}, we adopt a single scale parameter. The key is therefore to \textit{choose an effective scale parameter that enables the 2SNN to capture the significant scales.}

To this end, we choose the \textit{effective scale parameter as the geometric mean of all small parameters.} 
The idea is inspired by effective-medium theories for composite materials \cite{Ransome2021,ShengGeoEffect1982}, where multiple physical coefficients are reduced to a single coefficient via geometric averaging. This provides a simple yet balanced way to represent a unified scale when several small parameters are present. Closely related ideas also arise in stabilized finite element methods \cite{Codina_stabilizedFEM_2000}, where stabilization parameters are chosen to approximate the influence of unresolved subgrid scales by balancing characteristic scales of diffusion, convection, and reaction.

Although conceptually simple, this extension substantially broadens its applicability to more complex and realistic multiscale systems.
Meanwhile, it preserves a streamlined architecture while enabling the network to capture sharp gradients and coupled fast-slow dynamics arising from multiple interacting scales. As a result, the proposed framework provides a robust and practical approach for tackling more complex multiscale problems without increasing architectural or training complexity. 

Below, we summarize our \emph{contributions} and the \emph{significance} of this work:

\begin{itemize}

\item We propose a \emph{streamlined} network-2SNN   \eqref{eq:nn-with-aux-general}, that augments the input with a scale-aware feature.
The design intrinsically accommodates multiple scales and facilitates the network to predict \emph{both magnitude and location} of large gradients without requiring additional special treatments; the loss formulation and training procedure are in line with those of standard PINNs.

\item The proposed framework can handle dynamical systems with \emph{multiple small parameters} (including high-contrast regimes) by enriching the network with a \emph{single effective scale parameter}. Its efficacy is demonstrated through numerical experiments on a range of dynamical systems.

\item We provide \emph{theoretical justification} for the network design and a \emph{rationale} for selecting the geometric mean as an effective scale parameter to represent multiple small parameters in the model.

\item We extend the curriculum learning scheme (successive training strategy) in \cite{twoscaleNN_2024} for effective and stable training of models with multiple small parameters.

\end{itemize}
%
%
To further clarify the novelty of the present work, we include the following Table~\ref{tab:compare_current_prior} to summarize the key distinctions and their significance across different aspects compared with prior work \cite{twoscaleNN_2024}.

\begin{table}[!htb]
\centering
\begin{adjustbox}{width=0.94\textwidth}
\begin{tabular}{p{2.5cm} p{4.5cm} p{4.5cm} p{5.3cm}}
\hline
\textbf{Aspect} & \textbf{Prior Work} & \textbf{This Work} & \textbf{Why It Matters} \\
\hline
\raggedright Problem class
& \raggedright Scalar problems with a single small parameter $\epsilon$
& \raggedright Coupled dynamical systems with multiple small parameters $\epsilon_i$
& \raggedright Extends applicability to more realistic multiscale systems
\tabularnewline \noalign{\vskip 4pt}

\raggedright Scale handling
& \raggedright Single $\epsilon$ directly embedded
& \raggedright Effective scale parameter (geometric mean of $\epsilon_i$)
& \raggedright Enables unified treatment of multiple scales
\tabularnewline \noalign{\vskip 4pt}

\raggedright Architecture
& \raggedright Scale-augmented feature for scalar problems
& \raggedright Scale-augmented feature tailored for coupled systems
& \raggedright Preserves streamlined design while handling system complexity
\tabularnewline \noalign{\vskip 4pt}

\raggedright Theory
& \raggedright Limited
& \raggedright Justification of effective scale choice and network design rationale
& \raggedright Provides theoretical support for the approach
\tabularnewline \noalign{\vskip 4pt}

\raggedright Training
& \raggedright Standard training or curriculum learning for single small parameter problems
& \raggedright Curriculum learning for multiple and high-contrast small parameter regimes
& \raggedright Provides a robust curriculum learning paradigm for problems with multiple small parameters
\tabularnewline
\hline
\end{tabular}
\end{adjustbox}
\vskip 3pt
\caption{Comparison between prior work and the proposed approach.}
\label{tab:compare_current_prior}
\end{table}
%


\subsection{Related works}\label{sec:related_works}

Solutions to singularly perturbed and multiscale problems often exhibit rapid variations that are difficult for standard neural networks to approximate, largely due to spectral bias in gradient-based training. In the training, low-frequency components are learned first, while high-frequency features are captured slowly. Improving accuracy often requires increased architectural \cite{wang2021eigenvector,Jin23} or training \cite{adaptweightBPINN23,LaiMultistage2024} complexity. 
Existing approaches can be broadly categorized into \emph{architecture design} and \emph{training strategies}.

\textbf{Network architectures.}
Many approaches introduce spectral enrichment or scale-aware representations. Fourier-feature methods \cite{wang2021eigenvector} augment inputs with trigonometric bases to explicitly encode high-frequency components, though requiring preselection of frequency ranges. Related input-scaling methods \cite{LiuCaiXu20,LiXuZhang20,LiXuZhang23} adopt multiple scaled inputs to represent different frequencies.
Other common approaches are based on singularly perturbed theory, such as asymptotic-preserving neural networks \cite{Jin23,LuWangWu23,BerLP22}, 
boundary-layer PINNs \cite{BLPINN2023}, singular-layer PINNs \cite{slPINN_EABE2025,SingLayerPINN2026}, fast-slow neural networks \cite{FSNN2025}, and slow invariant manifold methods \cite{SlowInvariant_Singursys2024, SlowInvariant_FastSlowODE2024}, which involve multiscale decompositions, analytical corrections, special transformations, or invariance properties.

Despite these advances, handling multiple small parameters remains challenging. Most methods rely on explicit scale decomposition or problem-specific modifications. In contrast, we introduce a single effective parameter to capture the aggregated multiscale effect while maintaining a simple architecture.

\textbf{Training strategies.}
Complementary approaches address multiscale behavior via training. Adaptive sampling \cite{multiadaptsamplePOF26,adaptsampleCWA24} and adaptive weighting \cite{anagnostopoulos2023residualbased,adaptweightBPINN23,DCPINNEuler23} focus training on regions with large residuals or gradients.
More generally, two-stage or progressive strategies first learn a coarse approximation and then refine. Examples include gradient boosting \cite{fang2023ensemble}, multilevel methods \cite{multilvNN23}, multi-stage training \cite{LaiMultistage2024}, and curriculum learning \cite{bengio2009curriculum,soviany2022curriculum}. Curriculum learning is particularly flexible and can be combined with lightweight sampling strategies \cite{munzer2022curriculum}, using prior information from earlier training stages \cite{twoscaleNN_2024}.
However, standard curriculum learning typically assumes a single parameter. Extensions to multiple parameters include using an effective parameter (as in this work) or hierarchical curricula \cite{HierCurriculum21,PAPINNspb25} that progressively incorporate different scales. 

 
\textbf{Paper organization.}
The remainder of this paper is organized as follows. In Section~\ref{sec:method}, we state the problems of interest and introduce the two-scale neural network architecture and formulation, and the selection of the effective parameter. 
The theoretical justification of the network design and the rationale of the selection of the effective scale parameter are discussed in Section \ref{sec:theory_justify}. 
Numerical results in Section \ref{sec:numerical_examples} are then presented to validate accuracy, robustness under multiple/high-contrast small parameters, and the effectiveness of the curriculum learning strategy. A brief conclusion and discussion are presented in Section \ref{sec:conclusion}.

\section{Problems of interest and methodology}\label{sec:method}
We consider an abstract form of singular perturbation dynamics:
for a temporal domain $D$, 
\begin{subequations}\label{problem_general}
\begin{align}
\mathcal{L}_\epsilon \mathbf{u}=\mathbf{f},~\text{in}~ D,\\
\mathbf{u}=\mathbf{g} ~\text{on}~ \partial D,\label{general_bc}
\end{align}
\end{subequations}
where $\mathbf{f}(\tau)=({f}_1(\tau), f_2(\tau), \ldots, f_n(\tau))^{\top}$ and 
$\mathbf{g}=({g}_1, g_2 \ldots, g_n)^{\top}$. 
With this abstract form, we consider both first-order and second-order singularly perturbed dynamics.
The \emph{first-order} dynamical system has the following form 
\begin{subequations}\label{model:first_order_ODE}
\begin{align}
    \frac{d\mathbf{u}_m}{d\tau} &= \mathbf{f}_m(\mathbf{u}(\tau), \boldsymbol{\epsilon}), \label{eq:ODE_nonsing}\\
    \Lambda\frac{d\mathbf{u}_{m+1:n}}{d\tau}  &= \mathbf{f}_{m+1:n}(\mathbf{u}(\tau), \boldsymbol{\epsilon}), 
\end{align}
together with the initial condition
\begin{equation}\label{eq:intial_cond}
\mathbf{u}(\tau_0)=\mathbf{u}_0,
\end{equation}
\end{subequations}
where $m<n$, 
$\boldsymbol{\epsilon}=(\epsilon_1,\epsilon_2,\cdots, \epsilon_{n-m})^\top$,
$\Lambda = \operatorname{diag}(\boldsymbol{\epsilon})$,
$\mathbf{u}=(u_1,u_2,\ldots,u_n)^{\top} $, 
$\mathbf{u}_m=(u_1,u_2,\ldots,u_{m})^{\top}$, 
$\mathbf{u}_{m+1:n}=(u_{m+1},u_{m+2},\ldots,u_{n})^{\top}$.
$\mathbf{f}_m$ and $\mathbf{f}_{m+1:n}$ are defined in the same manner. In this system, we may also consider the case of vanishing \eqref{eq:ODE_nonsing}.
%
In addition, we consider the  
\emph{second-order} singularly perturbed dynamical system 
\eqref{model:second_order_ODE} equipped with both initial and terminal conditions. 
For   
$D= (\tau_0,\tau_f)$, 
\begin{subequations} \label{model:second_order_ODE}
\begin{align}
    \frac{d^2}{d\tau^2}\mathbf{u}_m + A(\tau)\frac{d\mathbf{u}_m}{d\tau}&= \mathbf{f}_m(\mathbf{u}(\tau), \boldsymbol{\epsilon}), \label{eq:2ndODE_macro}\\
    \Lambda\frac{d^2\mathbf{u}_{m+1:n}}{d\tau^2} + B(\tau) \frac{d\mathbf{u}_{m+1:n}}{d\tau}  &= \mathbf{f}_{m+1:n}(\mathbf{u}(\tau), \boldsymbol{\epsilon}), 
\end{align}
together with the initial and terminal conditions
\begin{equation}\label{eq:initial_terminal_cond}
\mathbf{u}(\tau_0)=\mathbf{u}_0, \;\mathbf{u}(\tau_f)=\mathbf{u}_f,
\end{equation}
\end{subequations}
where $A(\tau), B(\tau)$ are time-dependent coefficient matrices. The case without \eqref{eq:2ndODE_macro} is also considered for this system.

\subsection{Two-scale neural network architecture for models with multiple small parameters }\hfill

We denote $D_c=(\tau_0, T)$ as the interior computational domain for the considered problems: $T$ is a modest positive number for the first-order dynamical system~\eqref{model:first_order_ODE}, and $T=\tau_f$ for the second-order dynamical system~\eqref{model:second_order_ODE}. 
Recall that the Two-Scale Neural Network (2SNN) \cite{twoscaleNN_2024} is constructed by augmenting the input of a feedforward neural network with the feature
$\boldsymbol{\phi} := \bigl(\epsilon^\gamma(\tau-\tau_c),\epsilon^\gamma\bigr)$, leading to the form
\begin{equation}\label{eq:nn-with-aux-general}
N(\tau, \boldsymbol{\phi})
:= N\bigl(\tau,\epsilon^\gamma(\tau-\tau_c),\epsilon^\gamma\bigr),
\qquad \gamma<0,
\end{equation}
where $\boldsymbol{\phi}$ encodes \emph{scale information} in the network input and consists of two scale-dependent components: a stretched coordinate relative to $\tau_c$ and a scale magnitude induced by the small parameter. 
The factor $\epsilon^\gamma$ controls how strongly the network magnifies
the fast transition near the reference location $\tau_c$.
For multiple parameters, 
$\epsilon$ should act as an effective parameter  aggregating 
scales and reflecting the combined influence of the multiple fast-slow dynamics
associated with different small parameters $\epsilon_i$.
In the current study, the network output $N \in \mathbb{R}^n$ approximates the solution $\mathbf{u}$, and
$\epsilon$ is chosen as an \emph{effective small parameter} representing multiple small parameters in the model, taken as their geometric mean:
\begin{equation}\label{eq:geo_mean_multi}
\epsilon=\sqrt[n]{\epsilon_1 \epsilon_2\cdots \epsilon_n}.
\end{equation}
Also, we take $\gamma=-1/2$ as in \cite{twoscaleNN_2024}.
The theoretical justification of the 2SNN design and the rationale of the selection of the effective parameter via \eqref{eq:geo_mean_multi} are provided in Section \ref{sec:theory_justify}.
%

%
%

Once an effective parameter is selected, the training strategy from our previous work \cite{twoscaleNN_2024} is applicable, in which we start with a relatively large parameter and gradually reduce it to the desired small-parameter regime. This approach, known as curriculum learning, is detailed in Section~\ref{sec:method_in_CL}.


\nopagebreak


\subsection{PINN formulation and curriculum learning}\hfill

Under the standard PINN framework (with which 2SNN is in line), the neural network solution $\mathbf{u}_\theta$ is obtained by minimizing a residual-based loss function. 
%
%
Similar to \cite{twoscaleNN_2024}, the discrete loss function is defined follows

\begin{align}
\mathcal{L}_{col}(\theta)
=
&\frac{1}{N_c}\sum_{i=1}^{N_c}\left|\mathbf{r}_{\theta}\left(t_r^i\right)\right|^2+\frac{\alpha}{N_b}\sum_{i=1}^{N_b}\left|\mathbf{u}_{\theta}\left(t_b^i\right)-\mathbf{g}(t_b^i)\right|^2 \label{eq:loss-general},
\end{align}
where $
 \mathbf{r}_\theta(\tau)= 
\mathcal{L}_{\epsilon}\mathbf{u}_\theta (\tau) - \mathbf{f}(\tau)$ and $\{{t}_r^i\}_{i=1}^{N_c}$ specify the collocation points in the interior computational domain $D_c$, $\{ t_b^i\}_{i=1}^{N_b}$ are collocation points on $\partial_D=\{\tau_0\}~\text{or}~\{\tau_0, \tau_f\}$, $\alpha \geq 1$ is adjustable weights.

\subsubsection{Curriculum learning}\label{sec:method_in_CL}\hfill

To improve prediction accuracy, we apply a curriculum learning scheme, referred to as a successive training strategy in \cite{twoscaleNN_2024},  to singularly perturbed models with \emph{multiple small parameters}.
This approach gradually optimizes neural network parameters (weights and biases) by starting from a model with larger scale parameters and progressively training toward the target parameter regime. The purpose is to avoid the strong sensitivity to the initial guess that often arises when the model is trained directly with very small parameters. The successive training procedure is summarized in Algorithm~\ref{alg:succesive-training}. A similar approach is used in \cite{PAPINNspb25} while no scale-aware feature
is used by alternating the two parameters considered therein.%

\begin{algorithm}[!htb]
\SetAlgoLined

\KwData{Training dataset, learning rates, effective parameter $\epsilon$ derived from the dynamical system model using the geometric mean in \eqref{eq:geo_mean_multi}, and other model parameters.}
\KwResult{Successively optimized weights and biases of the neural network.}

\textbf{Step 0. (Initialization and choice of starting scale)}\\
If $\epsilon$ is very small, select a moderate initial value $\epsilon_0$ (empirically, on the order of $10^{-1}$ or $10^{-2}$); otherwise set $\epsilon_0=\epsilon$.
Initialize all network weights and biases with Xavier initialization.

\medskip
\textbf{Step 1. (Training at the current parameter scale)}\\
Replace $\epsilon$ in the dynamical system by $\epsilon_0$, and train the 2SNN \eqref{eq:nn-with-aux-general} with the loss function in \eqref{eq:loss-general}.\\
Store the optimized weights and biases obtained with $\epsilon_0$.

\medskip
\textbf{Step 2. (Iterative continuation toward the target scale)}\\
Choose an integer $\ell>1$.\\
\While{$\epsilon_0 \ge \ell\,\epsilon$}{
    Update the effective parameter of the model with $\epsilon_0 \leftarrow \epsilon_0/\ell$;\\
    Initialize the network with the weights and biases obtained from the previous iteration;\\
    Choose the learning rates for the coming retraining round (empirically use smaller values when $\epsilon_0<10^{-2}$);\\
    Retrain the network following Step 1 with the updated $\epsilon_0$;\\
}
If $\epsilon_0 < \ell\,\epsilon$, set $\epsilon_0=\epsilon$, initialize the network with the trained weights and biases from the previous iteration, and implement Step 1 once more.

\medskip
\textbf{Step 3. (Termination)}\\
Stop when the prescribed maximal number of epochs is reached.

\caption{Successive training of 2SNN for dynamical systems with multiple small parameters \cite{twoscaleNN_2024}}
\label{alg:succesive-training}
\end{algorithm}

It should be noted that the successive training is not strictly necessary; rather, it serves as a practical and stable training paradigm when applying 2SNN to problems with very small parameters. 
{In fact, we provide numerical experiments demonstrating that training a complex model directly with small parameters can still achieve reasonable accuracy (see, e.g., Table~\ref{tab:error_FHN_1e-2_8_direct} and Figure~\ref{fig:results_FHN_1e-2_8_direct} in Example~\ref{exm:FNsys}).}

\subsection{Stretching and shifting parameters in 2SNN: mechanism and choice}\label{sec:justification}\hfill

The proposed 2SNN \eqref{eq:nn-with-aux-general} is expected to predict both the location and magnitude of sharp transitions in the solutions.
The augmented inputs $\epsilon^\gamma(\tau-\tau_c)$ and $\epsilon^\gamma$ are linearly combined through learnable weights and biases, allowing the network to generate multiple shifted and stretched responses rather than being tied to a single reference point $\tau_c$. During training, regions with large residuals dominate the loss, and thus these learned combinations are driven toward locations of sharp transitions. As a result, the network can learn both the locations and magnitudes of sharp solution transitions without knowing their positions in advance.

In this study, we choose $\gamma=-1/2$ in 2SNN for numerical implementations. The reason for this choice will be discussed in Section \ref{sub2:theoretic_intuition}.
Unless otherwise specified, $\tau_c$ is chosen as the midpoint of the computational temporal domain $D_c$ and serves solely as a fixed reference point in the stretched coordinate, without assuming any prior knowledge of the location of a large gradient (sharp transition) in the solution.

%


\section{Theoretical justification}\label{sec:theory_justify} 
We provide a theoretical rationale for the 2SNN design with \textit{a single effective parameter}, explaining its suitability for singularly perturbed systems with multiple small parameters.

\subsection{2SNN design: layer-width scaling and gradient amplification}\label{sub2:theoretic_intuition}

The architecture \eqref{eq:nn-with-aux-general} is motivated by classical singular perturbation theory. Solutions typically exhibit localized transitions (layers) over narrow intervals in $\tau$, where $O(1)$ variation occurs.

For a single small parameter $\epsilon$, the layer width is commonly $O(\epsilon^\beta)$ with $\beta>0$. Introducing the stretched coordinate
\[
\eta = \frac{\tau - \tau^*}{\epsilon^{\beta}},
\]
maps this region to an $O(1)$ domain. Writing $\mathbf{u}(\tau;\epsilon)=\mathbf{U}(\eta;\epsilon)$ and applying the chain rule gives
\begin{equation}\label{eq:uexact_deri_comp}
\partial_{\tau} \mathbf{u} = O(\epsilon^{-\beta}), 
\qquad
\partial_{\tau\tau} \mathbf{u} = O(\epsilon^{-2\beta}),
\end{equation}
which characterizes the large gradients. Typical values are $\beta=1$ or $\beta=\tfrac{1}{2}$, e.g. in \cite{ Linss2010LayerAdaptedMeshes,OMalley1991Spb}. 

To align this scaling with the network, define $\xi=\epsilon^{\gamma}(\tau-\tau_c)$ and write
\[
\mathbf{u}_{\theta}(\tau)=N(\tau,\xi,\epsilon^\gamma).
\]
A Taylor expansion in $\xi$ yields
\begin{align}\label{eq:TSE_NN_comp}
\mathbf{u}_{\theta}(\tau)
=
N(\tau,0,\epsilon^{\gamma})
+
\epsilon^{\gamma}(\tau-\tau_c)\partial_{\xi}N
+
\frac{1}{2}\epsilon^{2\gamma}(\tau-\tau_c)^2\partial_{\xi\xi}N
+O\!\big(\epsilon^{3\gamma}(\tau-\tau_c)^3\big).
\end{align}
We adopt the following structural assumption.

\begin{assu}[Single fast-scale structure]\label{assu:single_fast_scale_comp}
All rapid variation in $\tau$ is induced through $\xi=\epsilon^{\gamma}(\tau-\tau_c)$, while the network remains otherwise slowly varying in $\tau$.
\end{assu}

Differentiating \eqref{eq:TSE_NN_comp} gives
\begin{equation}\label{eq:uNN_deri_comp}
\partial_\tau \mathbf{u}_{\theta} = O(\epsilon^{\gamma}),
\qquad
\partial_{\tau\tau} \mathbf{u}_{\theta} = O(\epsilon^{2\gamma}).
\end{equation}
Matching \eqref{eq:uNN_deri_comp} with \eqref{eq:uexact_deri_comp} yields $\gamma=-\beta$.
Since $\beta\in\{1/2,1\}$, we take $\gamma=-1/2$. This milder stretching avoids excessive localization while still resolving sharp transitions, providing a balance between resolution and robustness.

\textbf{Extension to multiple parameters.}
Although derived for a single $\epsilon$, the formulation extends naturally to multiple small parameters. Coupled systems often exhibit an aggregated multiscale behavior; the augmented input in \eqref{eq:nn-with-aux-general} adapts the network resolution to this overall sharpness. This motivates representing multiple parameters by a single effective scale.

\subsection{Rationale for an effective parameter: balanced multiscale aggregation}\label{subsec:ODE_system_multi}

An effective parameter $\epsilon$ aggregates multiple $\epsilon_i$ to balance stability and resolution. It is not intended to resolve each scale individually, but to capture their combined effect.

If $\epsilon$ is too small, then $\epsilon^\gamma$ becomes large, causing the stretched coordinate $\epsilon^\gamma(\tau-\tau_c)$ to grow rapidly and destabilize training.
If $\epsilon$ is too large such that $\epsilon^\gamma=O(1)$,  the stretching effect vanishes, leading to under-resolution of sharp transitions.

The geometric mean \eqref{eq:geo_mean_multi} provides a balanced choice: it incorporates all parameters without being dominated by extremes, maintaining an appropriate level of stretching. This enables robustness under parameter contrast while preserving the ability to resolve sharp features.

Therefore, 2SNN with the geometric-mean effective parameter is expected to capture multiscale behavior in coupled systems. This choice will be validated numerically in Section \ref{sec:numerical_examples}.

\section{Numerical Examples}\label{sec:numerical_examples}

\paragraph{\textbf{Default settings for training, architecture, 2SNN and its corresponding $\epsilon$.}} In the numerical experiments, we adopt the $\tanh$ activation function, and Adam optimizer. The collocation points are sampled from a uniform distribution over the computational domain. The MLP used in the 2SNN for the numerical tests has the architecture $(3,10,10,10,10,n)$, with $n$ the dimension of the dynamical system. {This network size is chosen as an empirical balance between expressiveness and training speed; other choices are also viable.} The default piecewise learning-rate scheduler is defined in Table \ref{tab:lr_pt}.
Without stating otherwise, the effective $\epsilon$ of the model is obtained with \eqref{eq:geo_mean_multi}. We also assume that the measure of the time domain is of order one; if not, a suitable variable transformation is applied to rescale it accordingly. As a result, we adopt the standardized time domain $[0,1]$ for the temporal variable $\tau$. Therefore, unless otherwise specified, the \emph{2SNN used in the numerical experiments} takes the form $N(\tau, (\tau-0.5)/\sqrt{\epsilon}, 1/\sqrt{\epsilon})$.


\begin{table}[!htb]
\centering
\begin{tabular}{c|ccccc}
\hline
Training steps 
& $\le 10{,}000$ 
& $10{,}000$--$30{,}000$ 
& $30{,}000$--$50{,}000$ 
& $50{,}000$--$70{,}000$ 
& $\ge 70{,}000$ \\
\hline
Learning rate 
& $1\times10^{-3}$ 
& $5\times10^{-3}$ 
& $1\times10^{-3}$ 
& $5\times10^{-4}$ 
& $1\times10^{-4}$ \\
\hline
\end{tabular}
\caption{Default piecewise learning-rate scheduler.}
\label{tab:lr_pt}
\end{table}

\paragraph{\textbf{Test problems and verification of the key claims.}} We test and discuss the performance of the proposed 2SNN on the following four examples: a first-order (Example \ref{exm:IVP}) and a second-order system (Example \ref{exm:BVP}) with a single small parameter, as well as more complex systems with multiple small parameters (Examples \ref{exm:Robertson} and \ref{exm:FNsys}). 

 Through these examples, we aim to substantiate the following key claims for the proposed 2SNN. 
 Claim A concerns the ability to handle coupling effects among components in single-parameter systems, which lays a necessary foundation for the multiple small-parameter regime. Claims B, C, and D concern robustness to parameter contrasts, the role of effective parameter selection, and a stable training paradigm for the multiple small-parameter regime, respectively.

\textbf{Claim A}: It effectively addresses the \emph{coupled} system and accurately captures sharp solution transitions introduced by a small parameter (Examples \ref{exm:IVP} and \ref{exm:BVP}).

\textbf{Claim B}: It delivers robust predictions even under \emph{high-contrast} small-parameter regimes (Example \ref{exm:Robertson}).

\textbf{Claim C}: The choice of the effective parameter plays a critical role (Example \ref{exm:Robertson}).

\textbf{Claim D}: Curriculum learning enhances training stability (Example \ref{exm:FNsys}).

\subsection{Dynamical systems  with a single small parameter}\label{subsec:system_w_single}\hfill

To illustrate \textbf{Claim A}, we consider a first-order  (Example \ref{exm:IVP}) and a second-order fast-slow  (Example \ref{exm:BVP}) dynamical system with a single small parameter. These serve as canonical models exemplifying the general formulations in \eqref{model:first_order_ODE} and \eqref{model:second_order_ODE}, which commonly arise in enzymatic reaction kinetics \cite{MichaelisMenten1913, Hunter2004Asymptotic}.

The reference solution for Example \ref{exm:IVP} is computed using the classical fourth-order Runge-Kutta finite difference method with a very small step size of $10^{-5}$, while that for Example \ref{exm:BVP} is given analytically. The results will be evaluated using graphical comparisons with reference solutions and absolute and relative errors.

\begin{figure}[h]
	\centering
	\subfigure[reference and NN solution of $u$]
      {\includegraphics[width=0.28\textwidth]
        {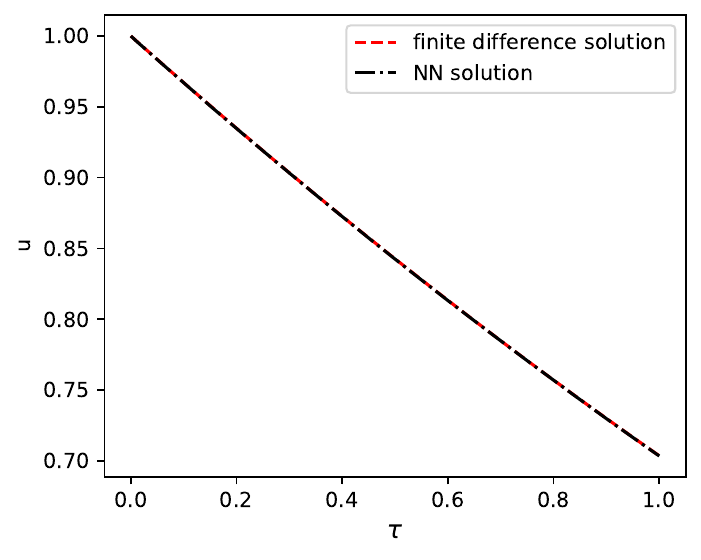}
        }
        \subfigure[absolute error of $u$]
        {\includegraphics[width=0.28\textwidth]{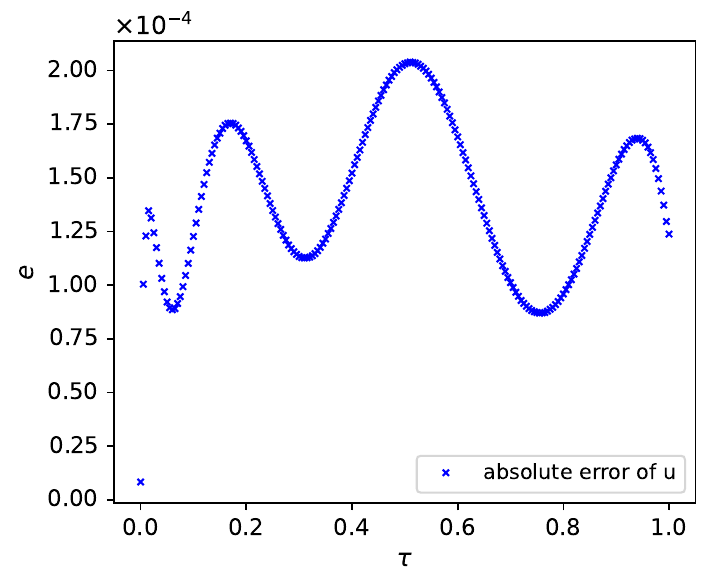}}
        \subfigure[relative error of $u$]
        {\includegraphics[width=0.28\textwidth]{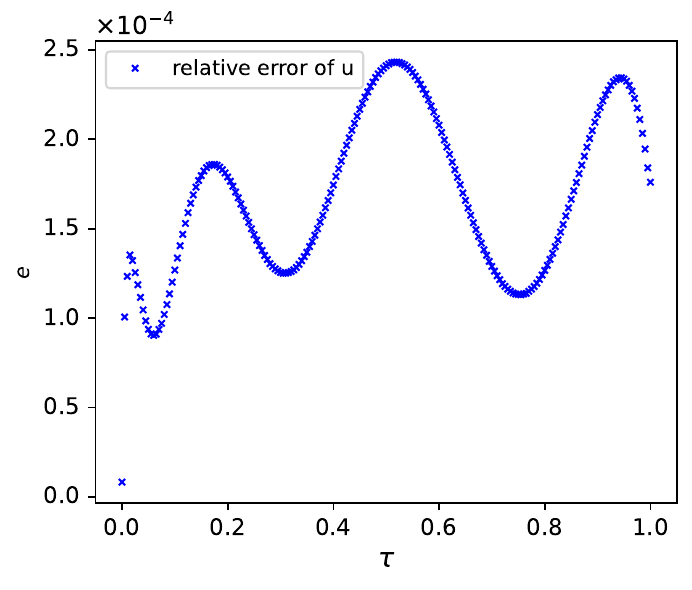}
        }

        \subfigure[reference and NN solution of $v$]
      {\includegraphics[width=0.28\textwidth]
        {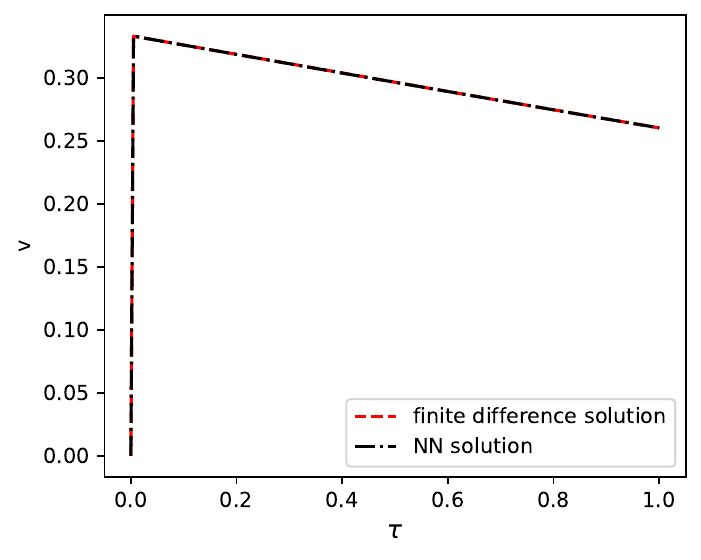}
        }
        \subfigure[absolute error of $v$]
        {\includegraphics[width=0.28\textwidth]{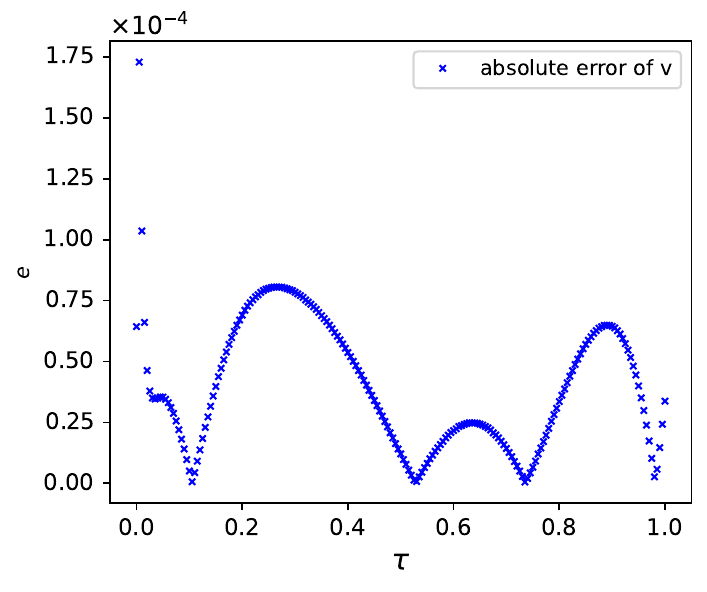}}
        \subfigure[relative error of $v$]
        {
        \includegraphics[width=0.28\textwidth]{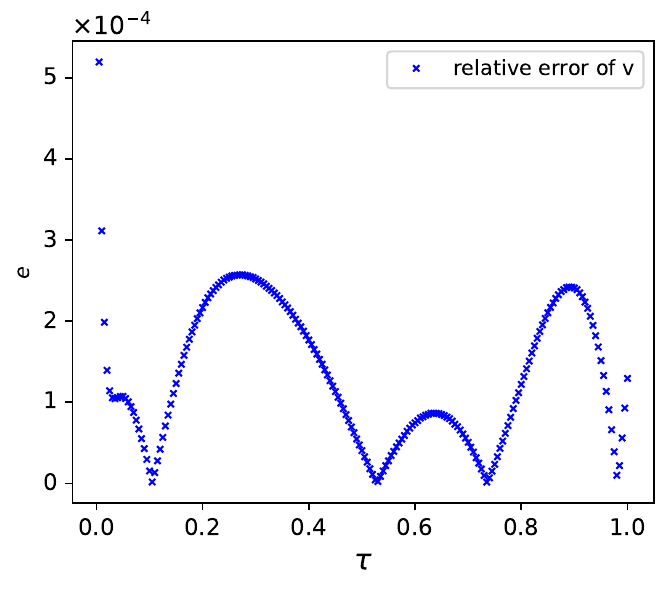}
        }
        \caption{Results for Example \ref{exm:IVP}  when $\epsilon=10^{-4}$
        { using 2SNN, with parameters specified in Table \ref{tab:exm1p1_para}  }.}
 \label{fig:results_IVP_1e_4}
        
\end{figure}

\begin{figure}[!htb]
	\centering
	\subfigure[reference and NN solution of $u$]
      {\includegraphics[width=0.29\textwidth]
        {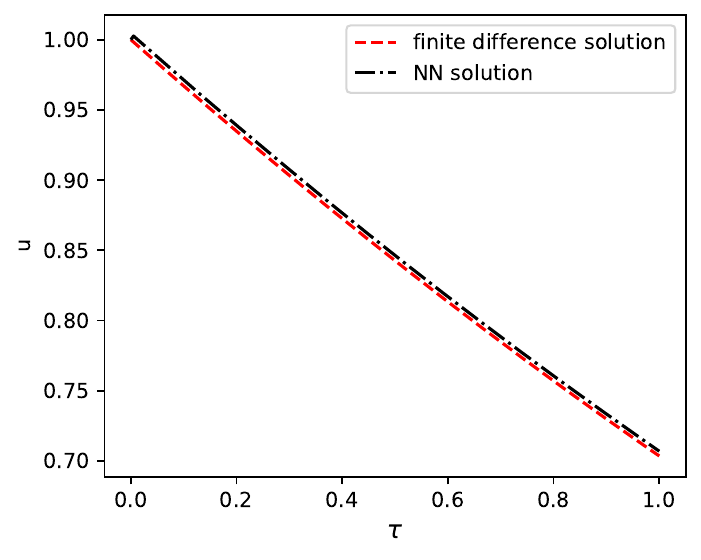}
        }
        \subfigure[absolute error of $u$]
        {\includegraphics[width=0.28\textwidth]{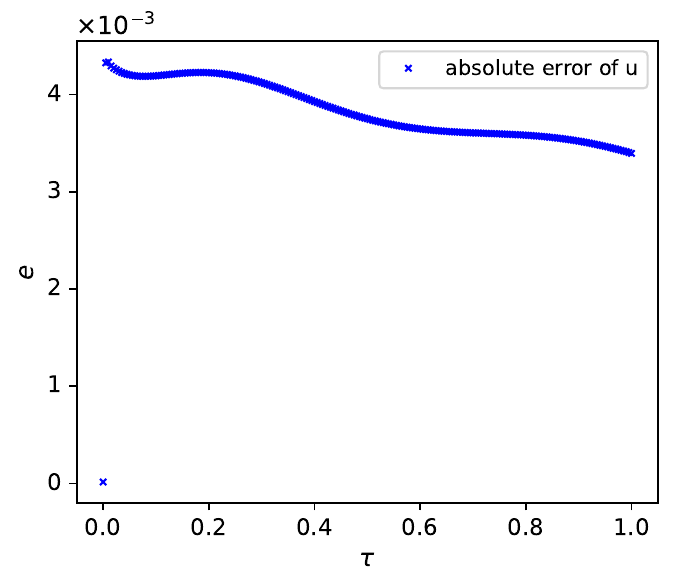}}
        \subfigure[relative error of $u$]
        {\includegraphics[width=0.28\textwidth]{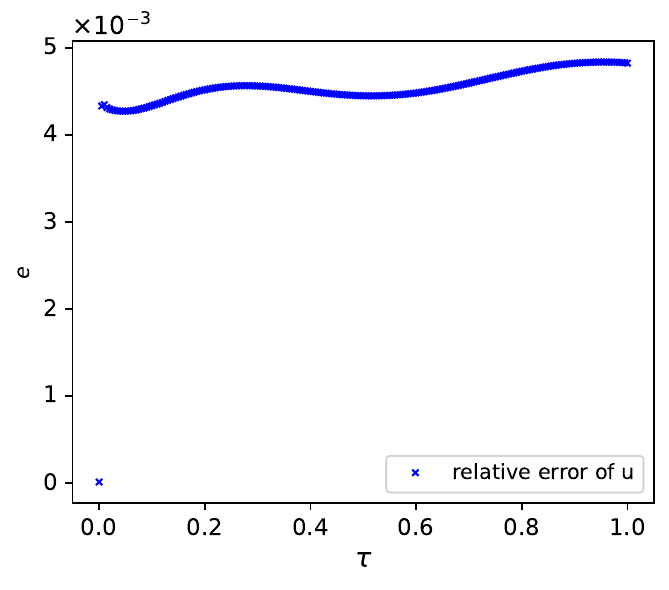}
        }

        \subfigure[reference and NN solution of $v$]
      {\includegraphics[width=0.29\textwidth]
        {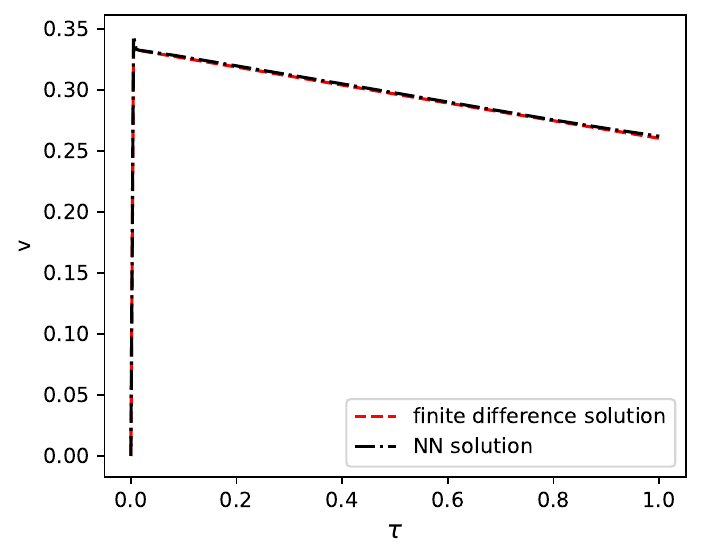}
        }
        \subfigure[absolute error of $v$]
        {\includegraphics[width=0.28\textwidth]{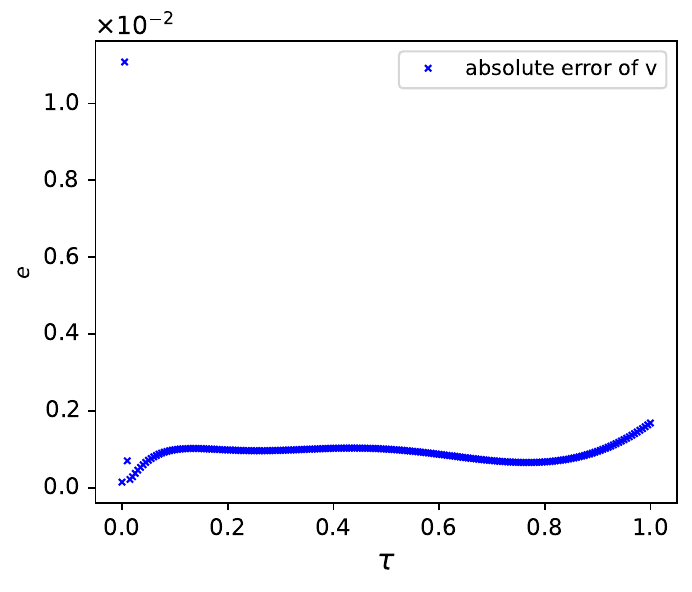}}
        \subfigure[relative error of $v$]
        {
        \includegraphics[width=0.28\textwidth]{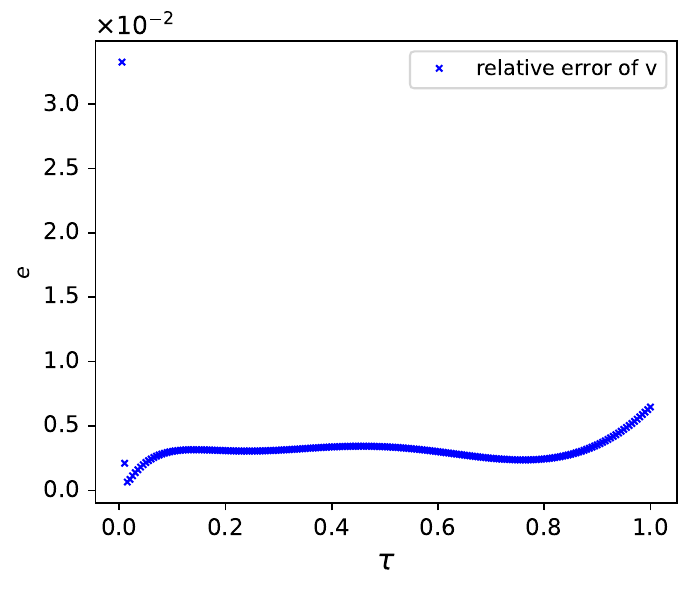}
        }
        \caption{Results for Example \ref{exm:IVP}  when $\epsilon=1.25\times 10^{-5}$
        { using 2SNN, with parameters specified in Table \ref{tab:exm1p1_para}  }.}
 \label{fig:results_IVP_1e_4_8}      
\end{figure}

\begin{figure}[!htb]
	\centering
     \subfigure[from $\epsilon=10^{-1}$ to $10^{-2}$.]{
        \includegraphics[width=0.26\textwidth]{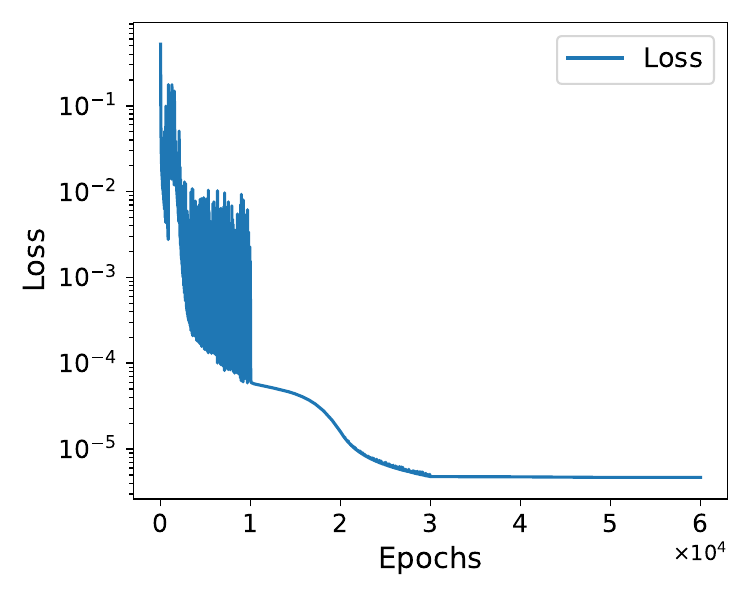}
        }
        \subfigure[$\epsilon=10^{-2}$ to $10^{-3}$]{
        \includegraphics[width=0.26\textwidth]{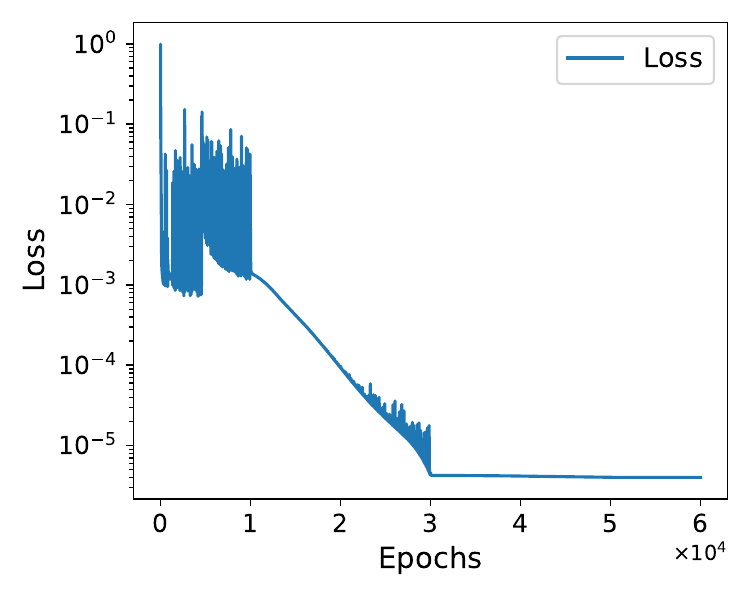}
        }
        \subfigure[$\epsilon=10^{-3}$ to $10^{-4}$]{
        \includegraphics[width=0.26\textwidth]{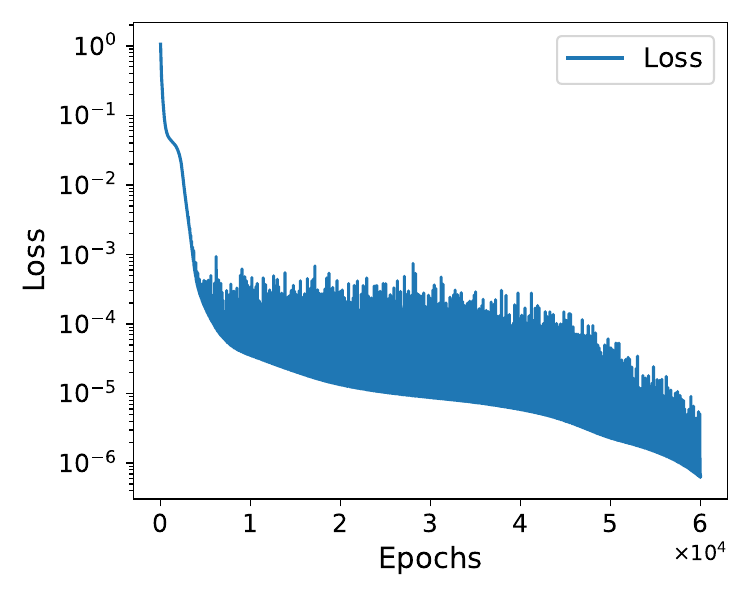}
        }
        
        \subfigure[$\epsilon=10^{-4}$ to $5\times10^{-5}$]{
        \includegraphics[width=0.26\textwidth]{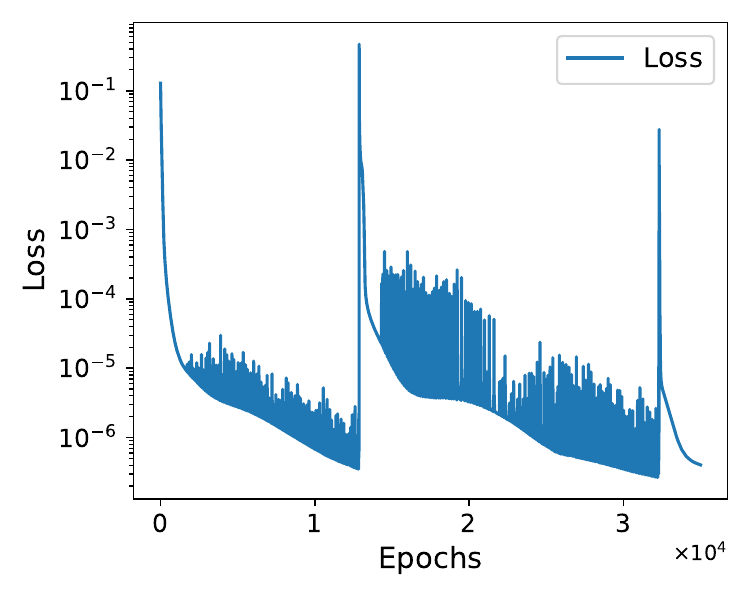}
        }
        \subfigure[$\epsilon=5\times10^{-5}$ to $2.5\times 10^{-5}$]{
        \includegraphics[width=0.26\textwidth]{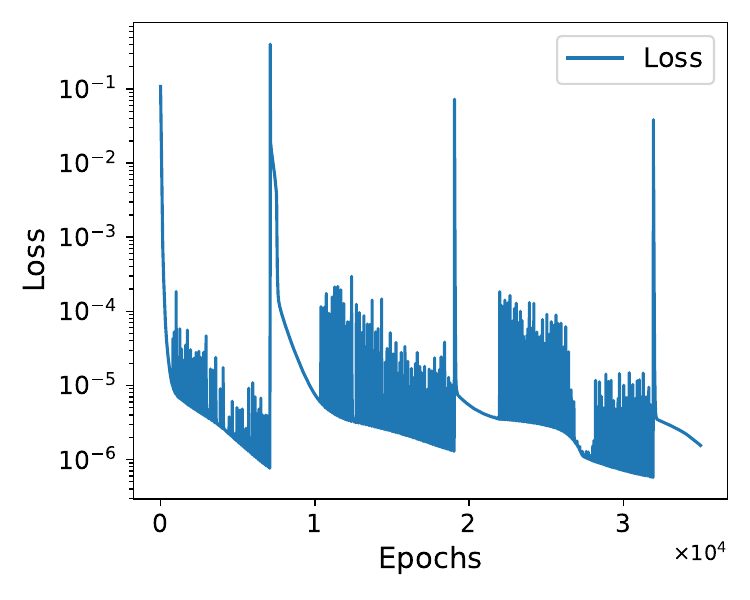}
        }
        \subfigure[$\epsilon=2.5\times10^{-5}$ to $1.25\times10^{-5}$]{
        \includegraphics[width=0.26\textwidth]{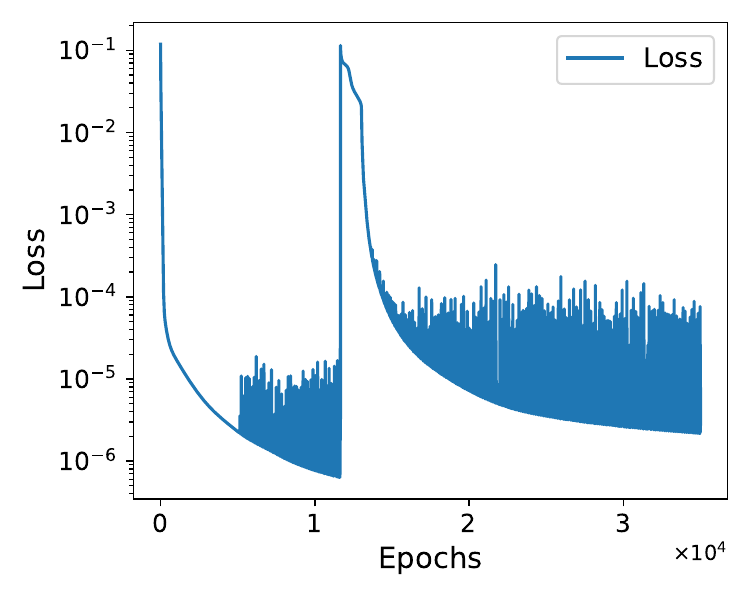}
        }
        \caption{ Loss history for Example \ref{exm:IVP} using successive training strategy with parameters in Table \ref{tab:exm1p1_para}.
}
		 \label{fig:IVP_MSNN_loss}
\end{figure}

\begin{exm}[First-order system: single small parameter] \label{exm:IVP}
\begin{subequations}\label{prob:IVP1_ori}
\begin{align}
\frac{ds}{dt} &= -k_1 \epsilon_0 s + (k_1 s + k_0) c, \\
\frac{dc}{dt} &= k_1 \epsilon_0 s - (k_1 s + k_0 + k_2) c, \\
s(0)& = s_0, \;c(0) = 0.
\end{align}
\end{subequations}
Introducing the following transforms:
\begin{align*}
u(\tau) &= \frac{s(t)}{s_0}, v(\tau) = \frac{c(t)}{\epsilon_0}, \tau = k_1 \epsilon_0 t,
\end{align*}
and let
$\lambda = \frac{k_2}{k_1 s_0}, k = \frac{k_0 + k_2}{k_1 s_0}, \epsilon = \frac{\epsilon_0}{s_0}$, then the problem \eqref{prob:IVP1_ori} can be rewritten as: to solve $(u,v)$ from
\begin{subequations}
\begin{align}
\frac{du}{d\tau} &= -u + (u + k - \lambda) v, \\
\epsilon\frac{dv}{d\tau} &=  u - (u + k) v, \\
u(0) & = 1, v(0) = 0.
\end{align}
\end{subequations}
where $\lambda > 0$ and $k > \lambda$.
We use $\lambda=1,k=2$ for the numerical test.
\end{exm}
We utilize the successive training in Algorithm \ref{alg:succesive-training} with the initial $\epsilon_0=10^{-1}$, and training toward $\epsilon=1.25\times 10^{-5}$ with the parameters and intermediate $\epsilon$ values specified in  Table \ref{tab:exm1p1_para}.
We collect the results for $\epsilon = 10^{-4}$ and $\epsilon = 1.25 \times 10^{-5}$ in
Figures~\ref{fig:results_IVP_1e_4} and \ref{fig:results_IVP_1e_4_8}, respectively. 
In both cases, the 2SNN provides accurate predictions.
When $\epsilon = 10^{-4}$, the relative errors of both $u$ and $v$ are on the order
of $10^{-4}$. When $\epsilon = 1.25 \times 10^{-5}$, although the prediction of $u$
shows minor discrepancy in Figure \ref{fig:results_IVP_1e_4_8} (a), its relative error remains on the order of $10^{-3}$ as shown in Figure \ref{fig:results_IVP_1e_4_8} (c).
The steep initial layer (nearly a vertical line) of $v$ is captured very well, and the relative error of
$v$ remains on the order of $10^{-2}$. The loss history of the successive training process is provided in Figure \ref{fig:IVP_MSNN_loss}. 

\begin{table}[!htb]\label{tab:IVP_progress}
\centering
\begin{tabular}{c c c c c c}
\hline
$\epsilon$ & $10^{-1}$ (initial $\epsilon_0)$    & $10^{-2}$ & $10^{-3}$ & $10^{-4}$ & $\boldsymbol{\epsilon}_{s}$ \\
\hline
$\alpha$ & 1 & 1 & 1 & 1 & 1 \\
$N_c$ & 300 & 300 & 450 & 450 & 450 \\
LR & P-S & P-S & P-S & $10^{-4}$ & $10^{-4}$ \\
iterations & $3\times 10^4$ & $6\times 10^4$ & $6\times 10^4$ & $6\times 10^4$ & $3.5\times 10^4$ each \\
\hline
\end{tabular}
\caption{ Parameters in the loss function \eqref{eq:loss-general} and hyper-parameters of the successive training for Example \ref{exm:IVP}. Here, $\boldsymbol{\epsilon}_s=10^{-5}\times [5, 2.5, 1.25]$, LR stands for the learning rate, P-S is the piecewise constant scheduler in Table \ref{tab:lr_pt}.}
\label{tab:exm1p1_para}
\end{table}

\begin{figure}[!htb]
	\centering
	\subfigure[exact and NN solution of $u$]
      {\includegraphics[width=0.28\textwidth]{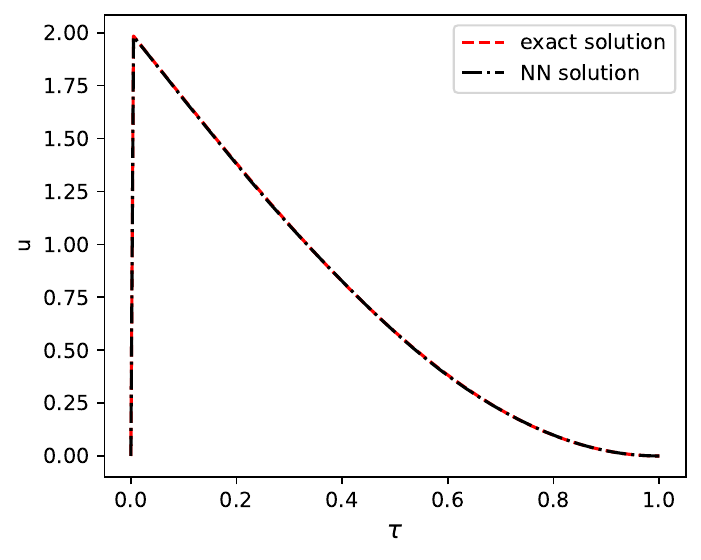}
        }
        \subfigure[absolute error of $u$]
        {\includegraphics[width=0.28\textwidth]{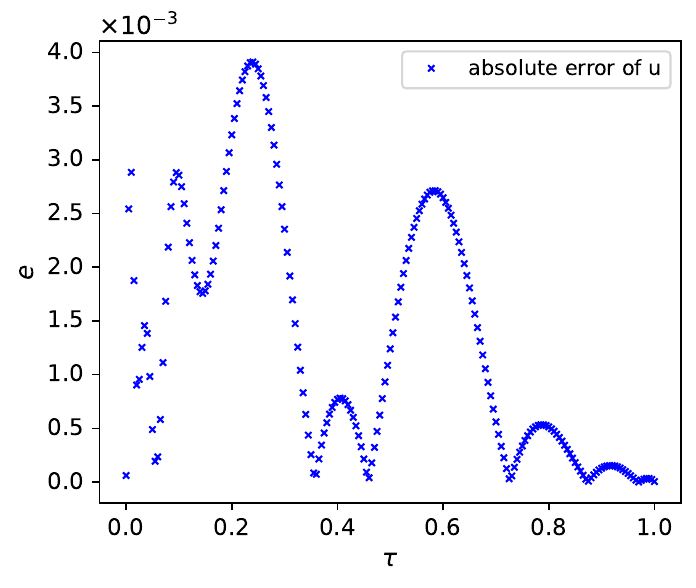}}
        \subfigure[relative error of $u$ around the initial layer]
        {\includegraphics[width=0.28\textwidth]{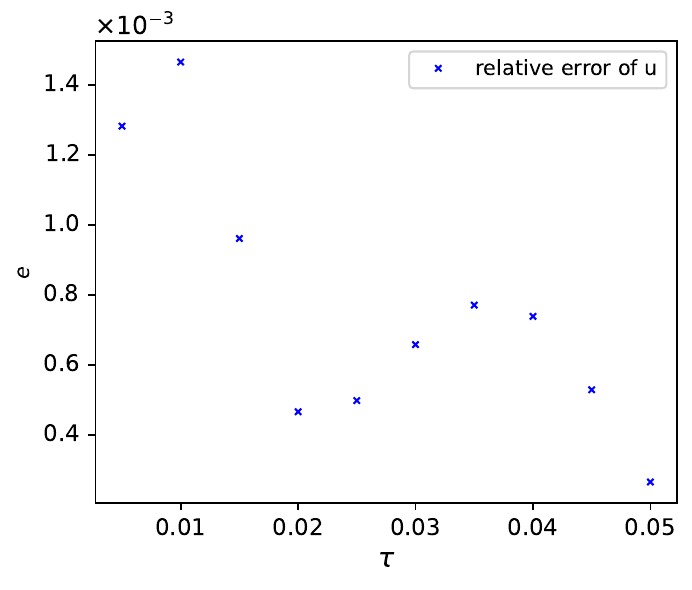}
        }

        \subfigure[exact and NN solution of $v$]
      {\includegraphics[width=0.28\textwidth]
        {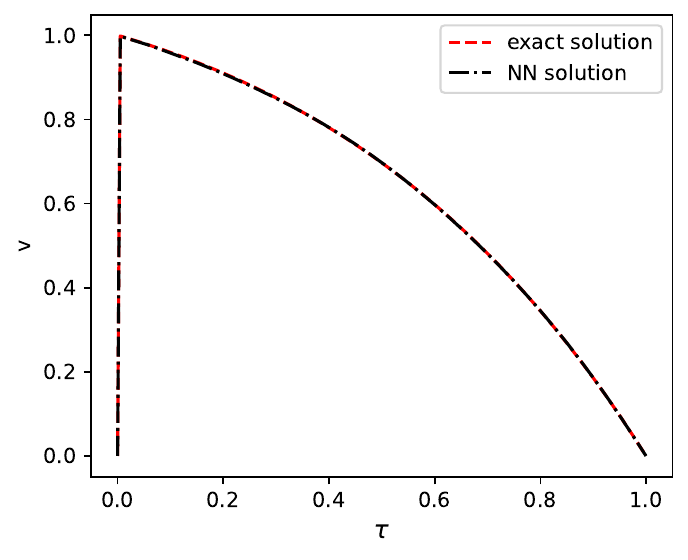}
        }
        \subfigure[absolute error of $v$]
        {\includegraphics[width=0.28\textwidth]{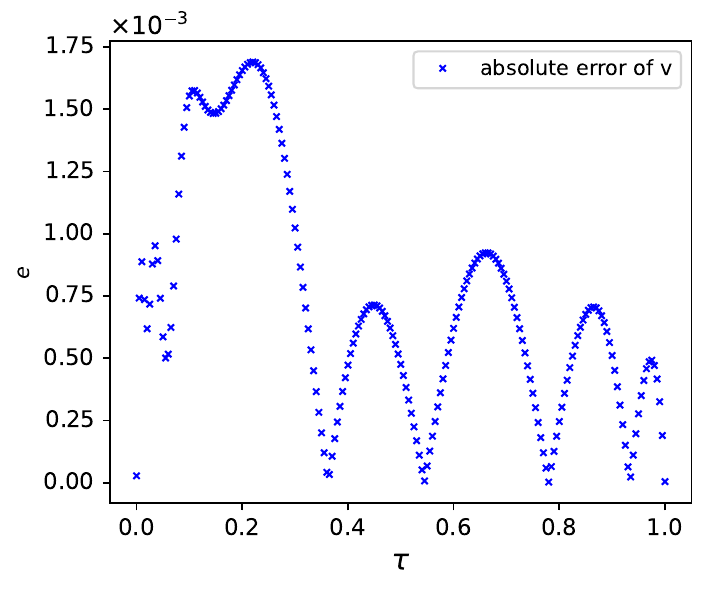}}
        \subfigure[relative error of $v$ around the initial layer]
        {
        \includegraphics[width=0.28\textwidth]{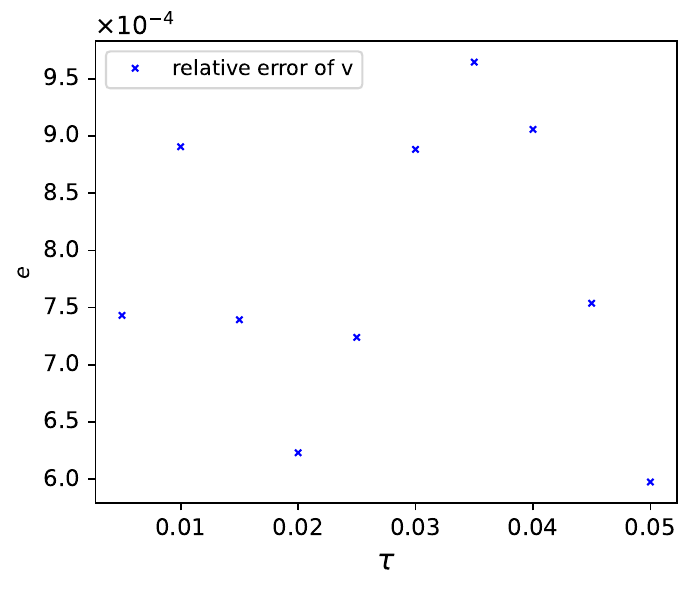}
        }
        \caption{Results for Example \ref{exm:BVP}  when $\epsilon=10^{-5}$ 
        {using 2SNN, with parameters specified in Table \ref{tab:exmbvp_para}}.}
 \label{fig:results_BVP_1e_5}
\end{figure}

\begin{figure}[!htb]
	\centering
     \subfigure[from $\epsilon=10^{-1}$ to $10^{-2}$.]{
        \includegraphics[width=0.22\textwidth]{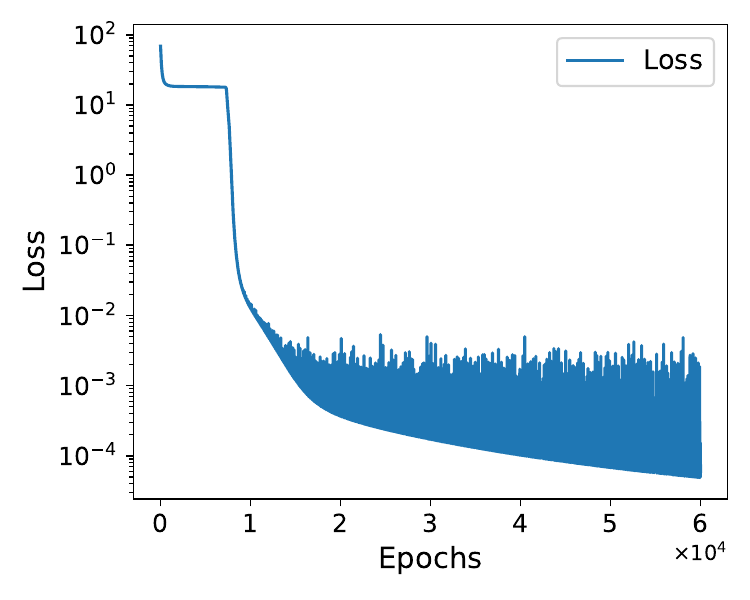}
        }
        \subfigure[$\epsilon=10^{-2}$ to $10^{-3}$]{
        \includegraphics[width=0.22\textwidth]{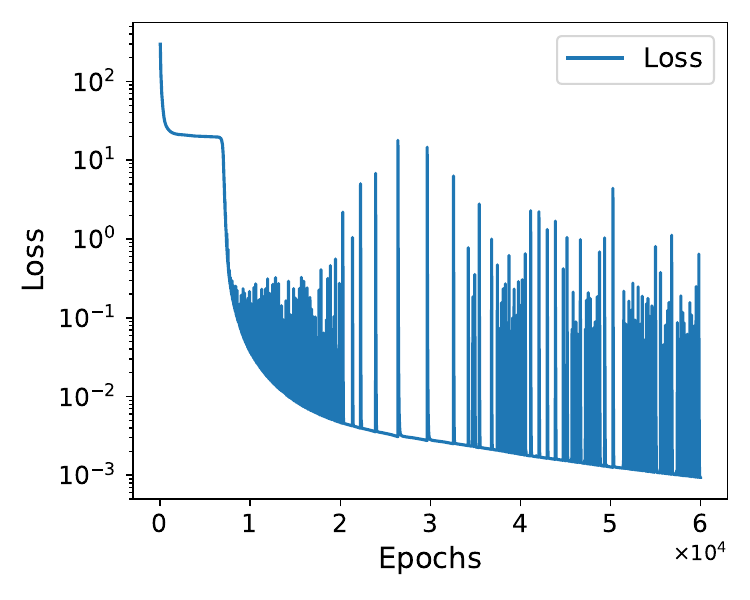}
        }
        \subfigure[$\epsilon=10^{-3}$ to $10^{-4}$]{
        \includegraphics[width=0.22\textwidth]{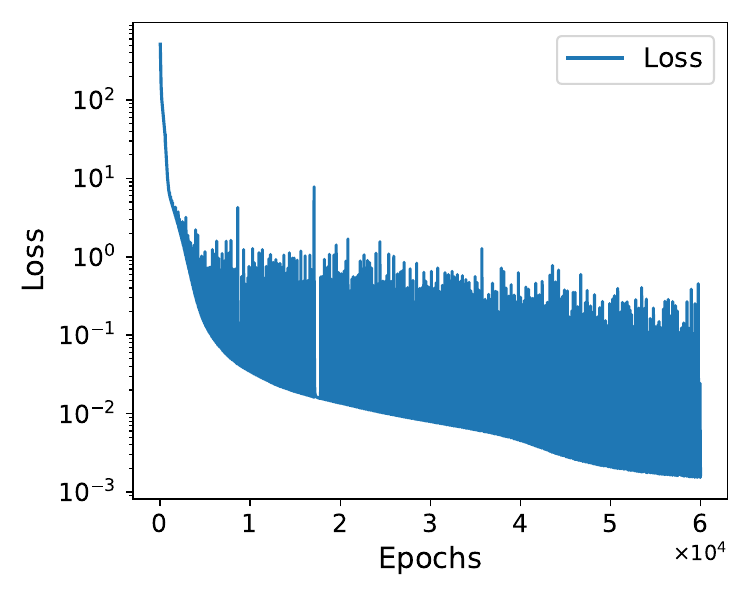}
        }
        \subfigure[$\epsilon=10^{-4}$ to $10^{-5}$]{
        \includegraphics[width=0.22\textwidth]{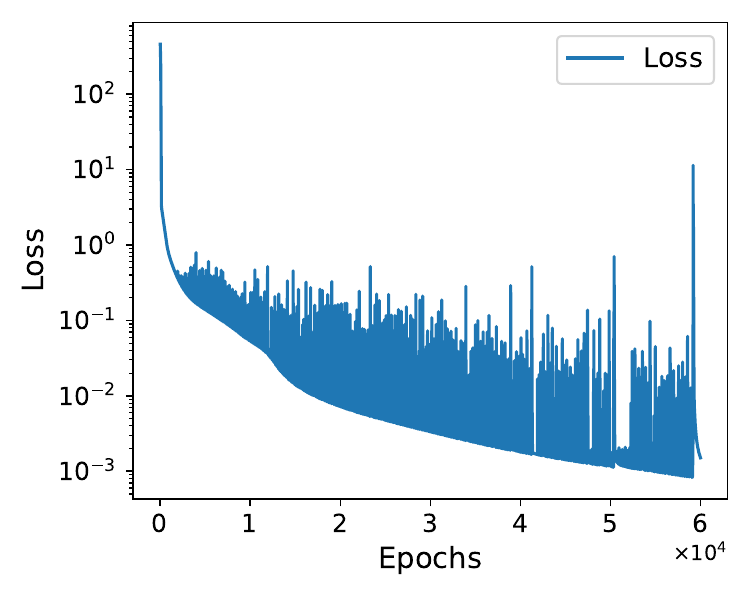}
        }
        \caption{ Loss history for Example \ref{exm:BVP} using successive training strategy with parameters in Table \ref{tab:exmbvp_para}.}
 \label{fig:BVP_MSNN_loss}
\end{figure}

\begin{exm}[Second-order system with initial and terminal conditions]\label{exm:BVP}
\begin{subequations}
\begin{align}
    \epsilon \frac{d^2u}{d\tau^2} + \frac{d u}{d\tau} - 2u + v = f_1, \\
    \epsilon \frac{d^2v}{d\tau^2} + 2\frac{d v}{d\tau} + u - 4v = f_2,
\end{align}
satisfying the initial and terminal conditions
\begin{align}
    u(0) = u(1) = v(0) = v(1) = 0,
\end{align}
\end{subequations}
where \( f_1(\tau) \), \( f_2(\tau) \), and the exact solution $u(\tau)$, $v(\tau)$ are given by
\begin{align*}
f_1(\tau) &=
\frac{4e^{-\frac{\tau}{\epsilon}} - \sin \left(\frac{\pi \tau}{2} \right) \pi^2 \epsilon^2 \left(1 - e^{-\frac{1}{\epsilon}}\right)}
{2\epsilon \left(-1 + e^{-\frac{1}{\epsilon}}\right)}
+ \frac{2e^{-\frac{\tau}{\epsilon}}}{\epsilon \left(1 - e^{-\frac{1}{\epsilon}}\right)}
- \pi \cos \left(\frac{\pi \tau}{2} \right)\\
&+ \frac{4 - 4e^{-\frac{\tau}{\epsilon}}}{-1 + e^{-\frac{1}{\epsilon}}}
+ 4 \sin \left(\frac{\pi \tau}{2} \right)
+ \frac{1 - e^{-\frac{2\tau}{\epsilon}}}{1 - e^{-\frac{2}{\epsilon}}}
- \tau e^{\tau-1},
\end{align*}

\begin{align*}
f_2(\tau) =
\frac{-4e^{-\frac{2\tau}{\epsilon}} - e^{\tau-1} \epsilon^2 \left(1 - e^{-\frac{2}{\epsilon}}\right)(\tau+2)}
{\epsilon \left(1 - e^{-\frac{2}{\epsilon}}\right)}
+ \frac{4e^{-\frac{2\tau}{\epsilon}}}{\epsilon \left(1 - e^{-\frac{2}{\epsilon}}\right)}\\
+ 2(\tau-1)e^{\tau-1}
+ \frac{2 - 2e^{-\frac{\tau}{\epsilon}}}{1 - e^{-\frac{1}{\epsilon}}}
- 2 \sin \left(\frac{\pi \tau}{2} \right)
+ \frac{4 - 4e^{-\frac{2\tau}{\epsilon}}}{-1 + e^{-\frac{2}{\epsilon}}},
\end{align*}
\[
u(\tau) = 2 \left(\frac{1 - e^{-\tau/\epsilon}}{1 - e^{-1/\epsilon}} - \sin \left(\frac{\pi \tau}{2} \right) \right),\,\,
v(\tau) = \frac{1 - e^{-2\tau/\epsilon}}{1 - e^{-2/\epsilon}} - \tau e^{\tau-1}.
\]
\end{exm}
For Example \ref{exm:BVP}, we utilize the successive training in Algorithm \ref{alg:succesive-training} with the initial $\epsilon_0=10^{-1}$, and training toward $\epsilon=10^{-5}$ with the parameters and intermediate $\epsilon$ values specified in  Table \ref{tab:exmbvp_para}.
%
We collect the results for $\epsilon = 10^{-5}$ in Figure~\ref{fig:results_BVP_1e_5}, which indicate that the 2SNN provides accurate predictions. The relative error of $u$ is on the order of $10^{-3}$, while that of $v$ is on the order of $10^{-4}$. For both components, the steep initial layers (nearly vertical lines) are captured very well. The loss history of the successive training process is summarized in Figure~\ref{fig:BVP_MSNN_loss}.

Examples \ref{exm:IVP} and \ref{exm:BVP} in Section \ref{subsec:system_w_single} have verified Claim A, showing that the proposed 2SNN effectively addresses coupled systems with a small parameter and accurately captures sharp solution transitions that exist in one or multiple components. 
This provides the necessary foundation for 2SNN to address more complex multi-parameter models.

\begin{table}[h]\label{tab:BVP_progress}
\centering
\begin{tabular}{c c c c}
\hline
$\epsilon$ & $10^{-1}$ (initial $\epsilon_0$) & $10^{-2}$ & $\boldsymbol{\epsilon_s}$ \\
\hline
$\alpha$ & 100 & 100 & 100 \\
$N_c$ & 300 & 300 & 450 \\
LR & P-S & $10^{-4}$ & $10^{-4}$ \\
iterations & $3\times 10^4$ & $6\times 10^4$ & $6\times 10^4$ \\
\hline
\end{tabular}
\caption{Parameters in the loss function \eqref{eq:loss-general} and hyper-parameters of the successive training for Example \ref{exm:BVP}. Here, $\boldsymbol{\epsilon}_s=[10^{-3}, 10^{-4}, 10^{-5}]$, LR stands for learning rate, P-S is the piecewise constant scheduler in Table \ref{tab:lr_pt}.}
\label{tab:exmbvp_para}
\end{table}

\subsection{Dynamical systems of multiple small parameters}\hfill

  For the dynamical system \eqref{model:first_order_ODE} with multiple small parameters, we study the Robertson model in chemical kinetics \cite{Robertson1966} in Example \ref{exm:Robertson} and the FitzHugh--Nagumo (FHN) model in neuoscience \cite{spbFNsys1981, ChimeraNNreview2019} in Example \ref{exm:FNsys}. Example \ref{exm:Robertson} primarily illustrates \textbf{Claims B} and \textbf{C}, while Example \ref{exm:FNsys} focuses on \textbf{Claims B} and \textbf{D}. 
%

For these two examples, the reference solutions are obtained using the Radau IIA implicit fifth-order Runge-Kutta method in Python, which is suited for systems exhibiting strong stiffness. 
The graphical comparison of the neural network and reference solutions will be presented, and errors will be measured in absolute and relative errors. Other metrics, $\Vert \cdot \Vert_\infty$ and $\Vert \cdot \Vert_{l^2}$, denoting the discrete $l^{\infty}$ and $l^2$ norms, will be used in the related tables.
$e_f$ represents the pointwise absolute error in a component function $f$.

\begin{exm}[Robertson model]\label{exm:Robertson}
\begin{subequations}
\begin{align}
\frac{dx}{dt} &= -k_1 x + k_3 y z, \\\label{eq:IVPmsp_x}
\frac{dy}{dt} &= k_1 x - k_2 y^2 - k_3 y z, \\
\frac{dz}{dt} &= k_2 y^2,\\\label{eq:IVPmsp_z}
x(0)&=1, y(0)=0, z(0)=0.
\end{align}
\end{subequations}
\end{exm}
The Robertson model describes multiple fast-slow chemical reactions \cite{Multiparaspb25}. Choosing an appropriate transform is essential for representing this complex singularly perturbed system. In this model, $k_1\ll k_3 \ll k_2$, to explicitly express a singularly-perturbed system, we follow \cite{Multiparaspb25} and let $\epsilon_1=k_1/k_2, \epsilon_2=k_3/k_2$. We also let $\tau = k_1 t$, which scales the long time interval corresponding to the slow rate $k_1$ to order~1. According to \cite{Multiparaspb25}, adding \eqref{eq:IVPmsp_x}-\eqref{eq:IVPmsp_z} and utilizing the initial condition, it follows $x+y+z=1$, therefore the model reduces to the following system in terms of $y,z$ only:
\begin{subequations}\label{model:Robertson}
\begin{align}
\frac{\epsilon_1}{\epsilon_2}\cdot\frac{dy}{d\tau}
&= \frac{\epsilon_1}{\epsilon_2}(1 - y - z)
   - \frac{1}{\epsilon_2} y^2
   -  y z, \label{eq:Rob1}
\\
\epsilon_1\frac{dz}{d\tau}
&=  y^2, \label{eq:Rob2} \\
y(0) & =0, z(0)=0.
\end{align}
\end{subequations}
We utilize the following configurations (i) and (ii)  to test the 2SNN for solving \eqref{model:Robertson}. \\
\textbf{Case (i)}: $k_1=4\times 10^{-2}, k_2= 10 , k_3=1$ $(\epsilon_1=4\times 10^{-3}, \epsilon_2=10^{-1})$.\\
\textbf{Case (ii)}: $k_1=4\times 10^{-2}, k_2=100,k_3=1$ $(\epsilon_1=4\times 10^{-4}, \epsilon_2=10^{-2})$.

In this multiple small-parameter Robertson system, the smallest model parameter is $\epsilon_1$ under the setups in Cases (i) and (ii), which can be as small as $4\times 10^{-4}$. The \emph{contrast} (ratio) of the small parameters in \eqref{eq:Rob1} and \eqref{eq:Rob2} is $1/\epsilon_2$, which could be as large as $100$ under the given setup. 
If we follow the default choice of effective $\epsilon$ as stated in \eqref{eq:geo_mean_multi}, then $\epsilon=\epsilon_1/\sqrt{\epsilon_2}$. In addition, when presenting the results, the $y$-component is scaled by $10^m$ ($m=1$ or $2$) to bring its magnitude close to 1, since $y$ is small; the corresponding absolute error is scaled accordingly to $10^m e_y$, making it more reflective of the prediction accuracy.

\begin{figure}[!htb]
	\centering
	\subfigure[vanilla PINN ($N(\tau)$)]
      {\includegraphics[width=0.29\textwidth]{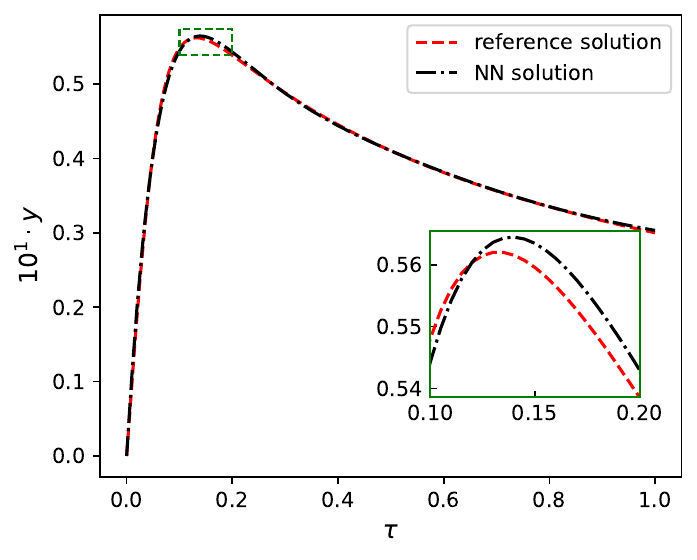}}
\subfigure[2SNN with $\epsilon=\epsilon_1/\sqrt{\epsilon_2}$]
        {\includegraphics[width=0.29\textwidth]{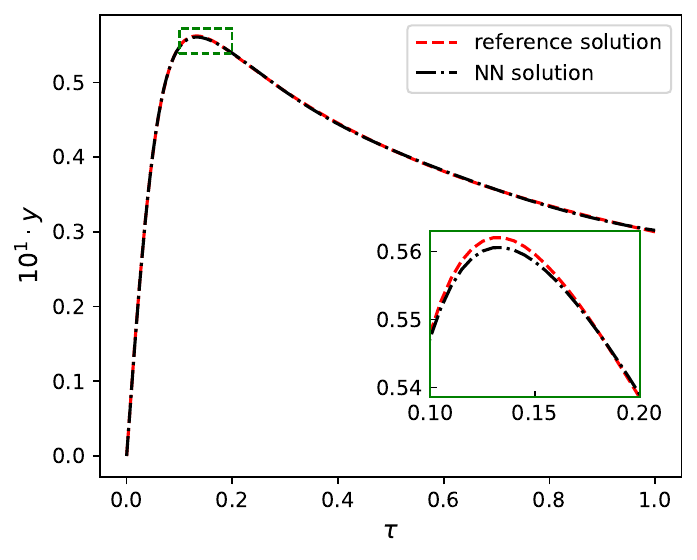}}
        \subfigure[2SNN with $\epsilon=\epsilon_1$]
        {\includegraphics[width=0.29\textwidth]{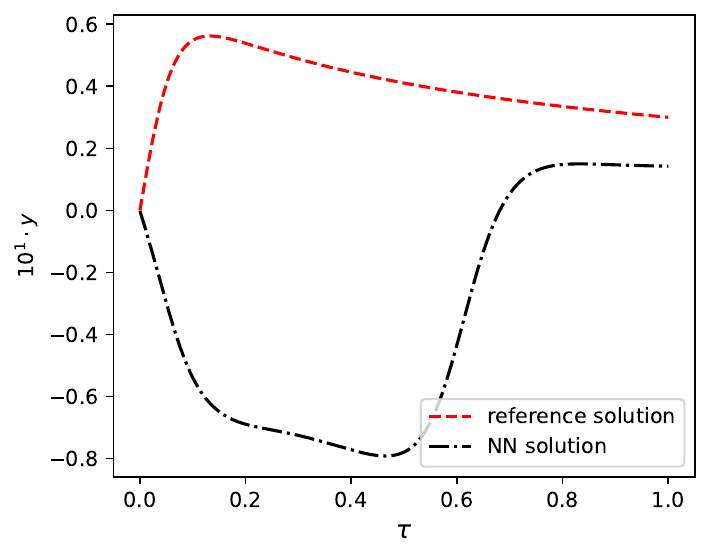}}
        \caption{Comparison of NN and reference solutions of $10y$ for Example \ref{exm:Robertson}  when $k_1=4\times 10^{-2}, k_2=10, k_3=1$, using $N(\tau)$ and $N(\tau, (\tau-0.5)/\sqrt{\epsilon}, 1/\sqrt{\epsilon})$ of different choices of $\epsilon$.}
 \label{fig:results_Rob_Casek2_10}
        
\end{figure}

For Case (i), we report the results in Figure~\ref{fig:results_Rob_Casek2_10} and Table~\ref{tab:error_Rob_k21050}. As shown in Figure~\ref{fig:results_Rob_Casek2_10} (b), the 2SNN using the default choice of effective $\epsilon$ ($\epsilon=\epsilon_1/\sqrt{\epsilon_2}$) accurately captures the reference solution. 
By comparison, the vanilla PINN (standard neural network), shown in Figure~\ref{fig:results_Rob_Casek2_10} (a), achieves a comparable accuracy to 2SNN. Nevertheless, the 2SNN remains more accurate, as evidenced by the zoom-in comparisons in Figure~\ref{fig:results_Rob_Casek2_10} (a) and (b), as well as by the third and fifth columns in the first two rows of Table~\ref{tab:error_Rob_k21050}.

We also explore another choice of scale parameter $\epsilon$ as a further validation of Claim C, where we take the smallest parameter of the model as the effective $\epsilon$, i.e., $\epsilon=\epsilon_1$. Under this setup, however, the 2SNN fails to capture the reference solution, as shown in Figure~\ref{fig:results_Rob_Casek2_10} (c) and in the fourth column of the first two rows of Table~\ref{tab:error_Rob_k21050}.
This degradation illustrates that an overly small effective scale parameter $\epsilon$ causes the stretched coordinate $(\tau-0.5)/\sqrt{\epsilon}$ in 2SNN to rapidly reach large values, thereby introducing instability, as mentioned in Section \ref{subsec:ODE_system_multi}.

\begin{figure}[!htb]
	\centering
	\subfigure[reference and NN solution of $10^2y$]
      {\includegraphics[width=0.28\textwidth]{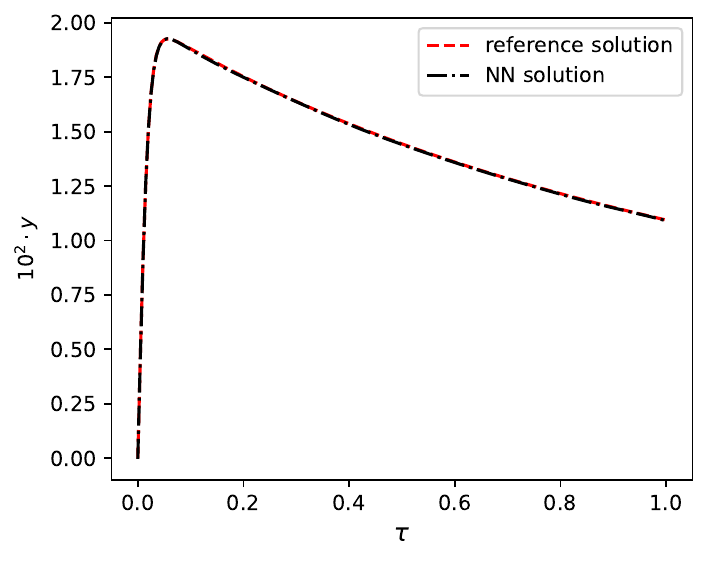}
        }
        \subfigure[absolute error of $10^2y$]
        {\includegraphics[width=0.28\textwidth]{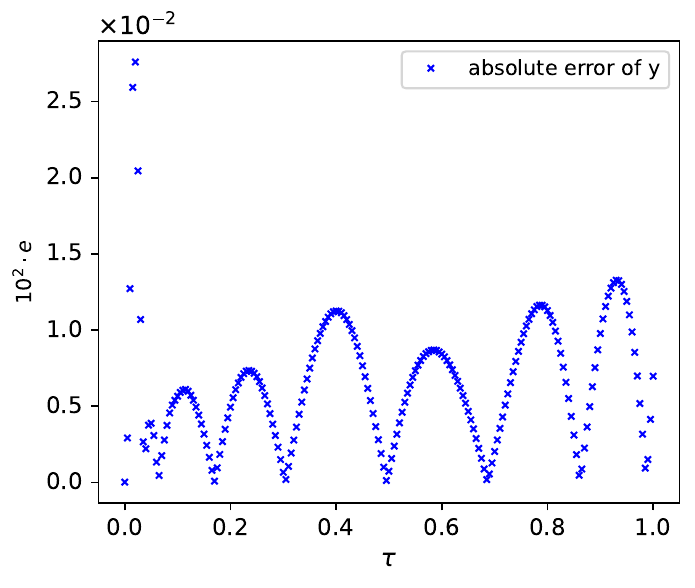}}
        \subfigure[relative error of $y$]
        {\includegraphics[width=0.28\textwidth]{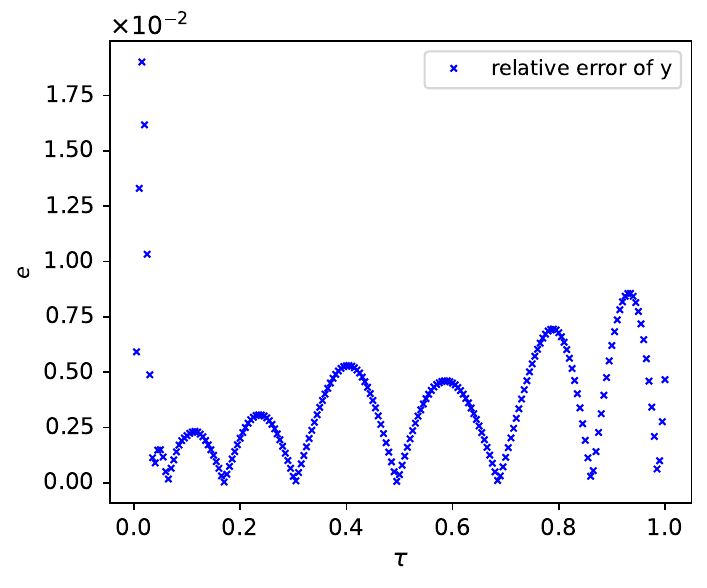}}

        \subfigure[reference and NN solution of $z$]
      {\includegraphics[width=0.28\textwidth]{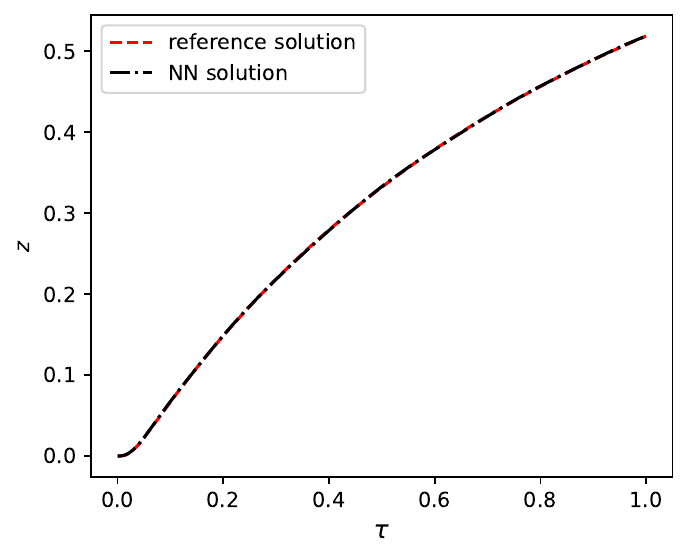}
        }
        \subfigure[absolute error of $z$]
        {\includegraphics[width=0.28\textwidth]{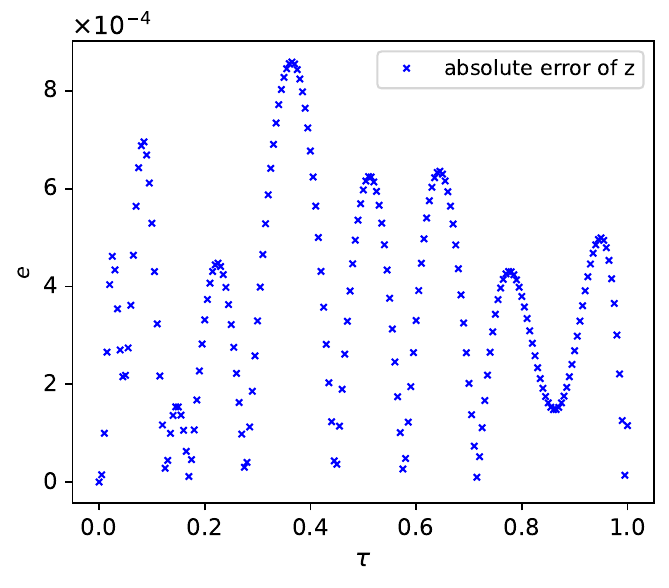}}
        \subfigure[relative error of $z$]
        {
        \includegraphics[width=0.28\textwidth]{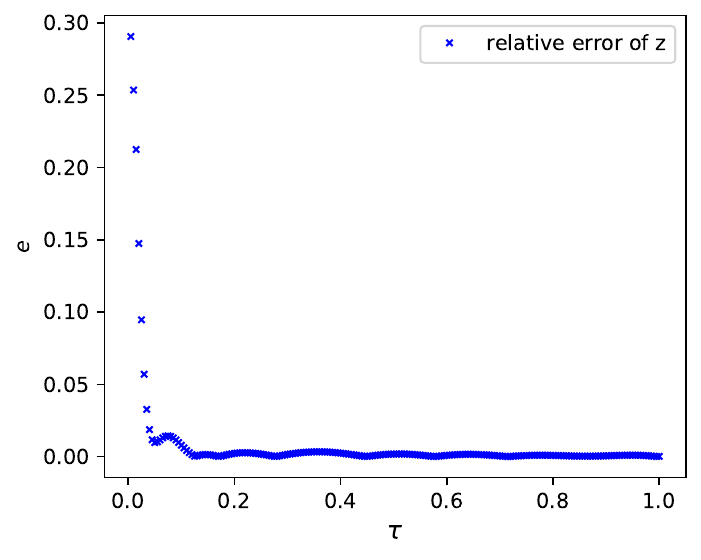}
        }
        \caption{Results for Example \ref{exm:Robertson}  when $k_1=4\times 10^{-2}, k_2=50, k_3=1$, using { $N(\tau, (\tau-0.5)/\sqrt{\epsilon}, 1/\sqrt{\epsilon})$}, with $\epsilon=\epsilon_1/\sqrt{\epsilon_2}$.}
 \label{fig:results_Rob_Casek2_50}
        
\end{figure}

%
We then implement the successive training in Algorithm \ref{alg:succesive-training} from the setup in Case (i) to solve \eqref{model:Robertson} toward the setup in Case (ii), with the parameter selection and intermediate $k_2$, as well as corresponding $\epsilon$ values listed in Table \ref{tab:Robt_successive}. Specifically, this successive training is regarding $k_2$, with $k_1$ and $k_3$ fixed.  
For the intermediate stage $k_2 = 50$, the results are presented in Figure~\ref{fig:results_Rob_Casek2_50} and Table~\ref{tab:error_Rob_k21050}. From Figure~\ref{fig:results_Rob_Casek2_50}, the 2SNN solution closely matches the reference solution, indicating an accurate prediction. The relative error of $z$ in Figure~\ref{fig:results_Rob_Casek2_50} (f) appears relatively large near $\tau=0$, which is because the solutions there are very close to zero and therefore highly sensitive to small numerical perturbations. 
%
We further compare the 2SNN and vanilla PINN results in Table~\ref{tab:error_Rob_k21050} under the same setup, where the vanilla PINN results are also computed with the successive training. The third and last columns of the table indicate that the 2SNN solution outperforms the vanilla PINN, achieving approximately one order of magnitude improvement in $10^{2} e_y$  in both the $l^\infty$ and $l^2$ norms. 

\begin{figure}[!htb]
	\centering
	\subfigure[reference and NN solution of $y$]
      {\includegraphics[width=0.28\textwidth]
        {figs/exm_robt/IVP_u_and_uNN_60to100.pdf}
        }
        \subfigure[absolute error of $10^2y$]
        {\includegraphics[width=0.28\textwidth]{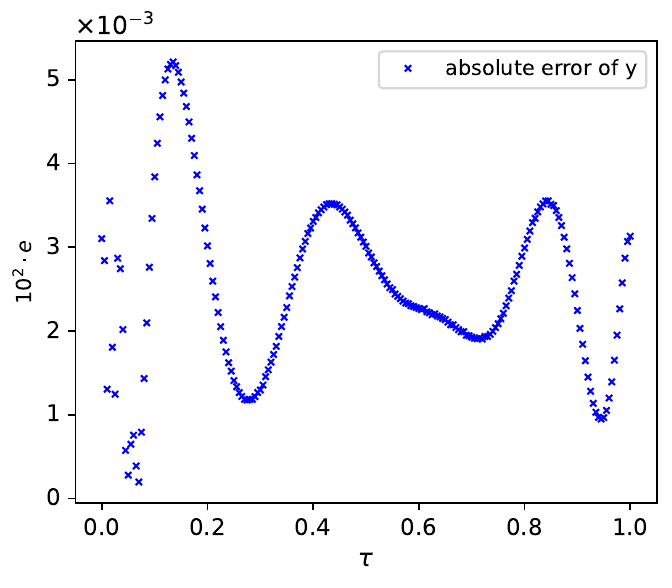}}
        \subfigure[relative error of $y$]
        {\includegraphics[width=0.28\textwidth]{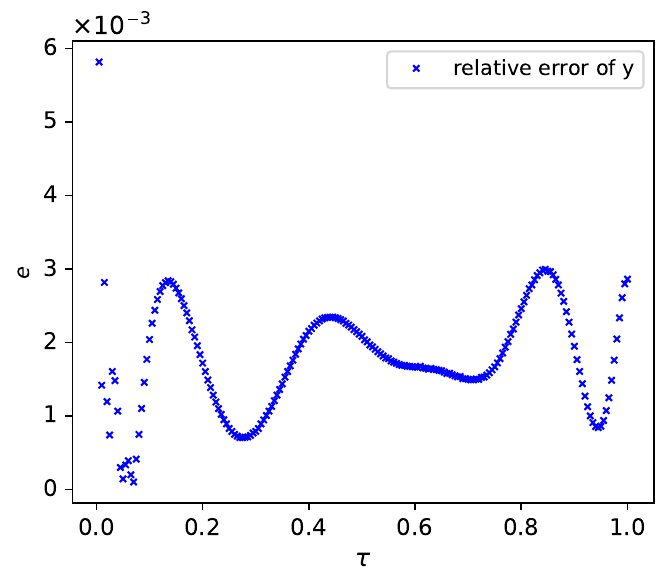}
        }

        \subfigure[reference and NN solution of $z$]
      {\includegraphics[width=0.28\textwidth]
        {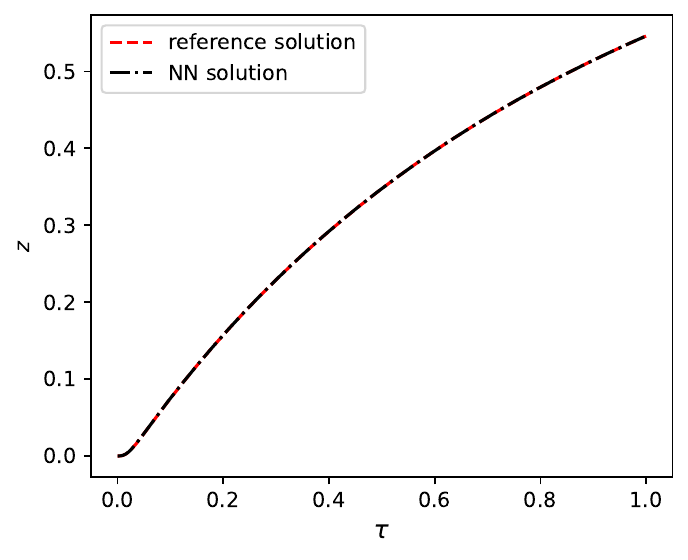}
        }
        \subfigure[absolute error of $z$]
        {\includegraphics[width=0.28\textwidth]{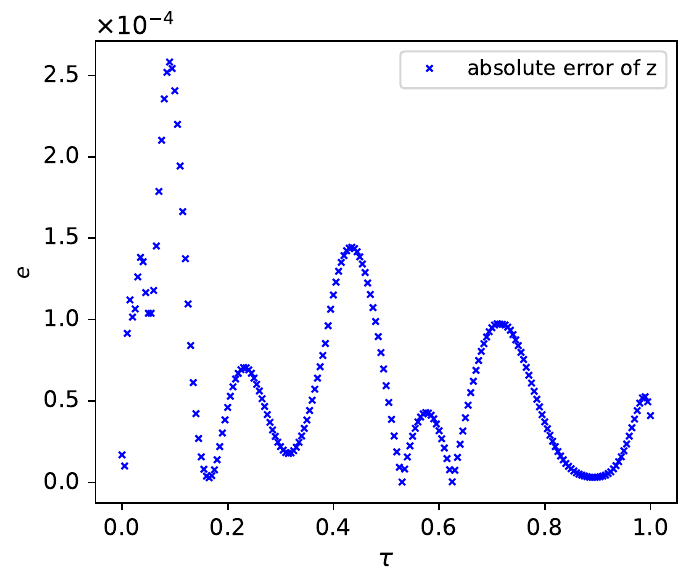}}
        \subfigure[relative error of $z$]
        {
        \includegraphics[width=0.28\textwidth]{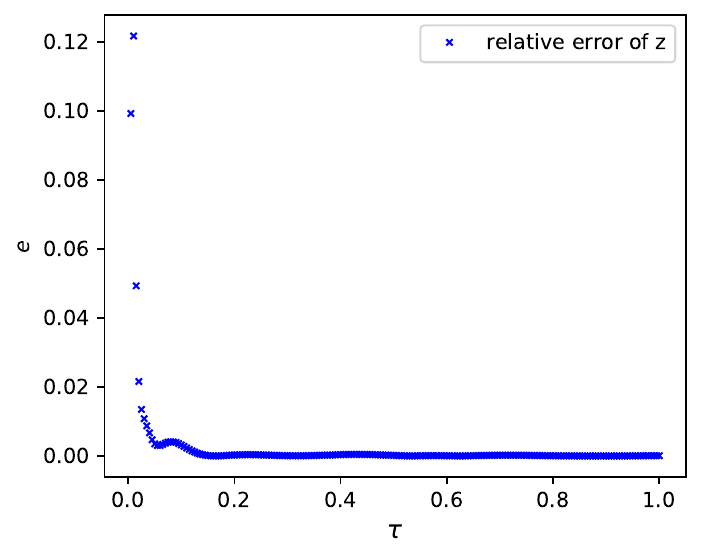}
        }
        \caption{Results for Example \ref{exm:Robertson}  when $k_1=4\times 10^{-2}, k_2=100, k_3=1$, using { $N(\tau, (\tau-0.5)/\sqrt{\epsilon}, 1/\sqrt{\epsilon})$}, with $\epsilon=\epsilon_1/\sqrt{\epsilon_2}$.}
 \label{fig:results_Rob_CaseII}
        
\end{figure}

The successive training finally arrives at the setup in Case~(ii), where $k_2=100$. The corresponding results are reported in Figure~\ref{fig:results_Rob_CaseII}, which again show excellent agreement between the 2SNN and reference solutions. As shown in Figure~\ref{fig:results_Rob_CaseII}    (c), the relative error of $y$ is on the order of $10^{-3}$. Figure~\ref{fig:results_Rob_CaseII}   (d)-(f) further demonstrate that 2SNN accurately captures the reference solution for $z$ as well. Similar to the previous discussion (for $k_2=50$), the larger relative error observed near $\tau=0$ in Figure~\ref{fig:results_Rob_CaseII}    (f) is because the solutions for $z$ there are very close to zero.

\begin{figure}[!htb]
	\centering
     \subfigure[$k_1=4\times 10^{-4}, k_3=1, k_2=50$ (trained from the same $k_1$ and $k_3$ but $k_2=10$)]{
        \includegraphics[width=0.28\textwidth]{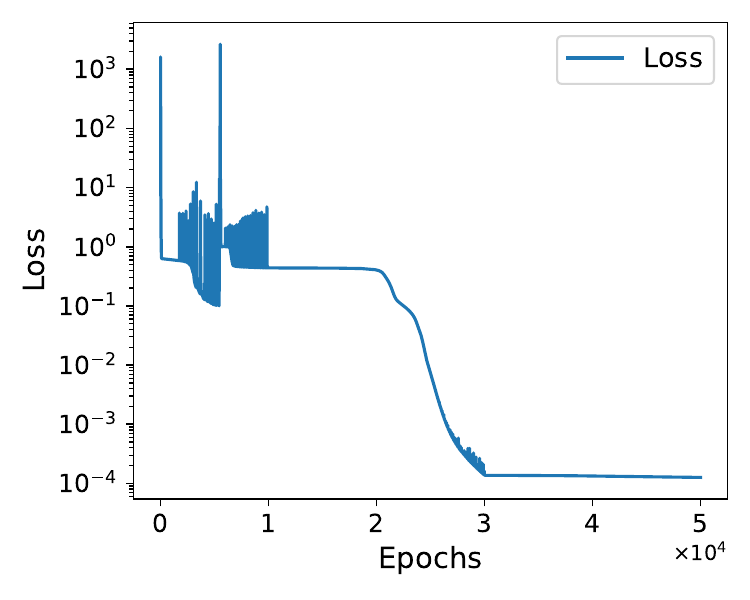}
        }
        \subfigure[$k_1=4\times 10^{-4}, k_3=1, k_2=100$ (trained from the same $k_1$ and $k_3$ but $k_2=60$)]{
        \includegraphics[width=0.28\textwidth]{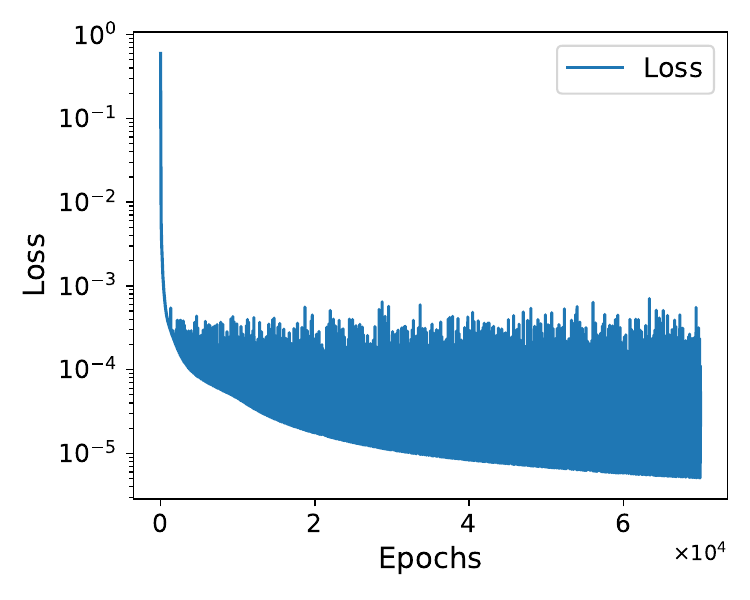}
        }

		\caption{ Loss history for Examples \ref{exm:Robertson}.
}
		 \label{fig:exmRobt_MSNN_loss}
\end{figure}

We also present the training loss history with respect to iterations (epochs) in Figure \ref{fig:exmRobt_MSNN_loss}, which records the training history from $k_2=10$ to $k_2=50$; and from $k_2=60$ to $k_2=100$.

\begin{table}[h]
\centering
 \begin{adjustbox}{width=0.59\textwidth}
\begin{tabular}{c c c c c}
\hline
$k_2$ & $10$ & $50$ & $60$ & $100$ \\
$\epsilon$ & $1.26\times 10^{-2}$ & $5.66\times 10^{-3}$ & $5.16\times 10^{-3}$ & $4\times 10^{-3}$ \\

\hline
$\alpha$   & $10^4$ & $10^4$ & $10^4$ & $10^4$ \\
$N_c$      & $300$  & $300$  & $300$  & $300$ \\
LR         & P-S    & P-S    & $10^{-4}$    & $10^{-4}$ \\
iterations & $5\times 10^4$ & $5\times 10^4$ & $5\times 10^4$ & $7\times 10^4$ \\
\hline
\end{tabular}
\end{adjustbox}
\caption{Parameters in the loss function \eqref{eq:loss-general} and hyper-parameters for the successive training with respect to $k_2$ for Example \ref{exm:Robertson} (with fixed $k_1=4\times 10^{-3}$ and $k_3=1$). Here, LR represents the learning rate,  P-S is the piecewise constant scheduler in Table \ref{tab:lr_pt}.}
\label{tab:Robt_successive}
\end{table}

\begin{table}[!htp]
    \centering
    \vspace{2mm}
    \renewcommand{\arraystretch}{1.1}
\begin{adjustbox}{width=0.63\textwidth}
    \begin{tabular}{c c|cc|c}
        \toprule
        & & \multicolumn{2}{c|}{{2SNN}} & Vanilla PINN \\
        \cmidrule(lr){3-4}
        \textbf{Case} & \textbf{Metric}
        & $\epsilon=\epsilon_1/\sqrt{\epsilon_2}$
        & $\epsilon=\epsilon_1$
        &  \\
        \midrule
        $k_2=10$ & ${\Vert 10 e_y \Vert}_{\infty}$
        & $1.97\times 10^{-3}$
        & $1.23$
        & $1.19\times 10^{-2}$ \\
         & ${\Vert 10 e_y \Vert}_{l^2}$
        & $7.63\times 10^{-4}$
        & $0.89$
        & $2.73\times 10^{-3}$ \\
        \midrule
        $k_2=50$ & ${\Vert 10^2 e_y \Vert}_{\infty}$
        & $2.77\times 10^{-2}$
        & not implemented
        & $2.86\times 10^{-1}$ \\
         & ${\Vert 10^2 e_y \Vert}_{l^2}$
        & $7.67\times 10^{-3}$
        & not implemented
        & $5.50\times 10^{-2}$ \\
        \bottomrule
    \end{tabular}
    \end{adjustbox}
    \vskip 3pt
\caption{errors of vanilla PINN and 2SNN solutions with different choices of effective $\epsilon$ for Example \ref{exm:Robertson}, with fixed $k_1= 4\times 10^{-3}, k_3=1.$}
\label{tab:error_Rob_k21050}
\end{table}

\begin{exm}[FitzHugh-Nagumo model]\label{exm:FNsys}
\begin{subequations}
\begin{align}
\epsilon_{1}\,\frac{d v}{d\tau}
&=
v - \dfrac{v ^3}{3} - z - w,\label{FHN_eq1}
\\
\epsilon_{2}\,\frac{d z}{d\tau}
&=
v - 0.5\,z,
\label{FHN_eq2}\\
\epsilon_{3}\,\frac{d w}{d\tau}
&=
v - w,
\end{align}
and the initial conditions  are given by
\begin{align}
v(0)=1.5, z(0)=0, w(0)=0.2. 
\end{align}
\end{subequations}
\end{exm}
In this example, we fix $\epsilon_1=10^{-1}$, $\epsilon_3=10^{-2}$, and alter $\epsilon_2$ gradually from  $\epsilon_2=10^{-2}$ to $\epsilon_2=2.5\times 10^{-4}$. 
Under this setup, the highest possible \emph{contrast} of the smallest parameter comes from \eqref{FHN_eq1} and \eqref{FHN_eq2} (i.e., $\epsilon_1/\epsilon_2$), which can be as large as $4\times 10^2.$
The specific intermediate values for $\epsilon_2$ during the successive training are recorded in Table \ref{tab:FNsys_para}, with the corresponding $\epsilon$ values. Under this setup, the smallest model parameter is $\epsilon_2$, which can be as small as $2.5\times 10^{-4}$.
When altering $\epsilon_2$, 
the effective  $\epsilon=\sqrt[3]{\epsilon_1 \epsilon_2 \epsilon_3}$ from \eqref{eq:geo_mean_multi} then   changes from $\epsilon=2.15\times 10^{-2}$ to $\epsilon=6.30\times 10^{-3}$ in the successive training process (Algorithm \ref{alg:succesive-training}). 
Besides implementing the successive training, we also directly solve for the case $\epsilon_1=10^{-1}$, $\epsilon_3=10^{-2}$, and $\epsilon_2=10^{-2}/8$ to further test the capacity of 2SNN. 

Under this parameter configuration, decreasing $\epsilon_2$ significantly sharpens the solution profiles of the state variable $y_2$, producing a pronounced initial layer.
A moderate value of $\epsilon_1$ suppresses oscillatory behavior, whereas a suitably small $\epsilon_3$ yields a tail that gradually flattens. Since our investigation focuses on singularly perturbed effects, this parameter configuration serves as a controlled setting that highlights the singular-perturbation behavior (i.e., initial layers) primarily governed by $\epsilon_2$, while suppressing oscillatory and sharp-tail behavior.

\begin{figure}[!htb]
	\centering
	\subfigure[reference and 2SNN solutions \quad of $v$]
      {\includegraphics[width=0.30\textwidth]{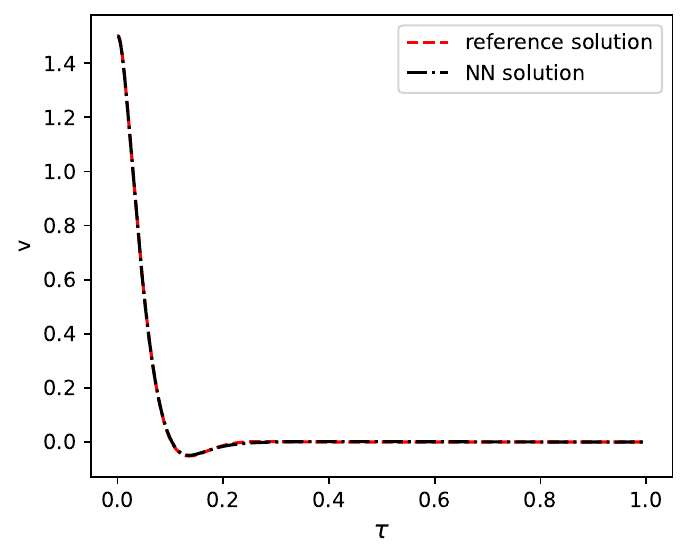}
        }
        \subfigure[reference and 2SNN solutions \quad of $z$]
        {\includegraphics[width=0.30\textwidth]{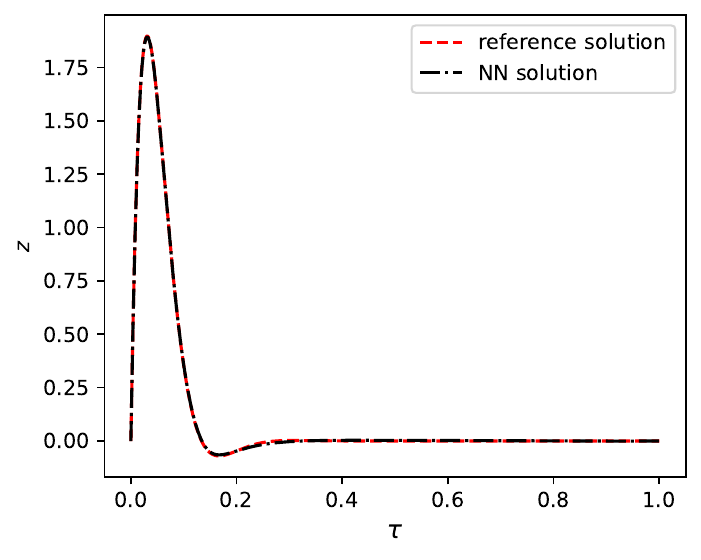}}
        \subfigure[reference and 2SNN solutions \quad of $w$]
        {\includegraphics[width=0.30\textwidth]{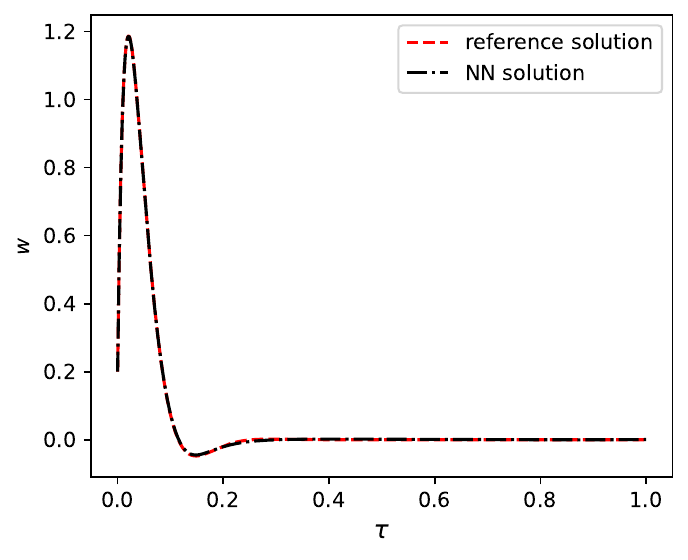}
        }

        \subfigure[absolute error of reference and 2SNN solutions of $v$]
      {\includegraphics[width=0.28\textwidth]{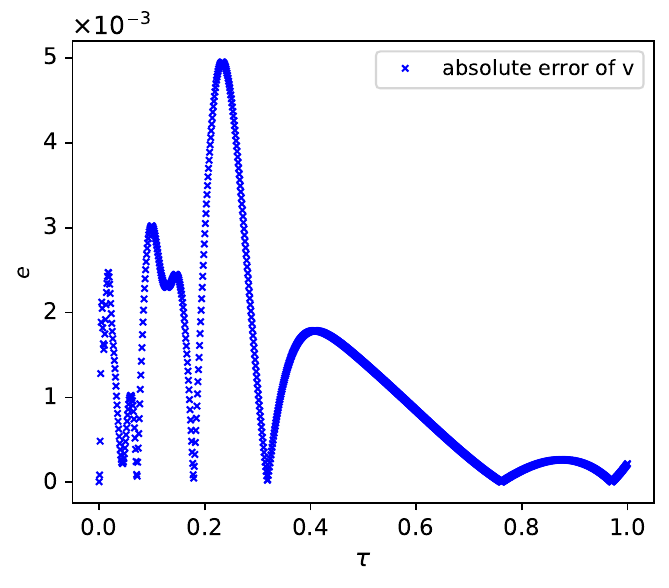}
       
        }
        \subfigure[absolute error of reference and 2SNN solutions of $z$]
        {\includegraphics[width=0.28\textwidth]{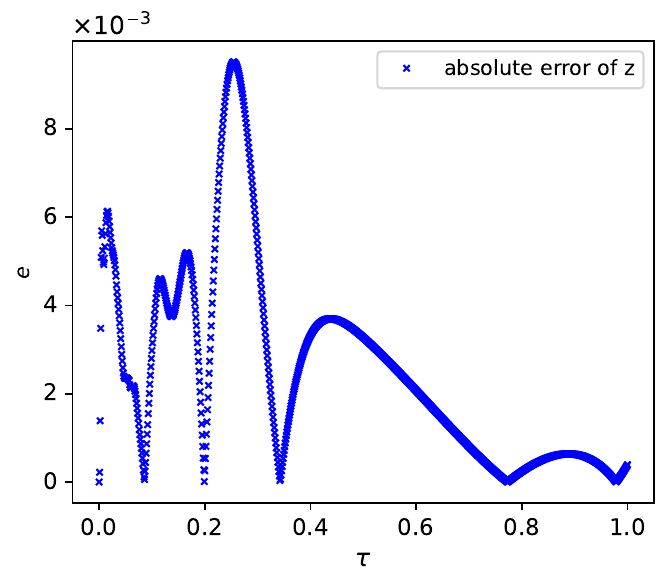}}
        \subfigure[absolute error of reference and 2SNN solutions of $w$]
        {
        \includegraphics[width=0.28\textwidth]{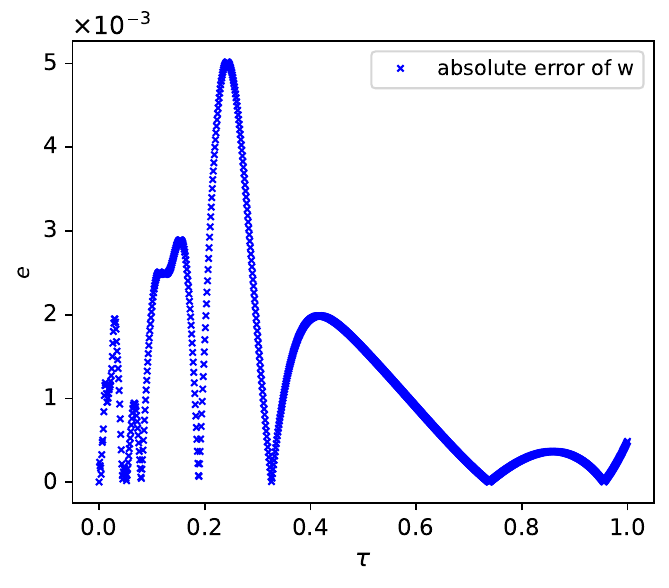}
        }
        \subfigure[absolute error of reference and vanilla PINN solutions of $v$]
      {\includegraphics[width=0.28\textwidth]{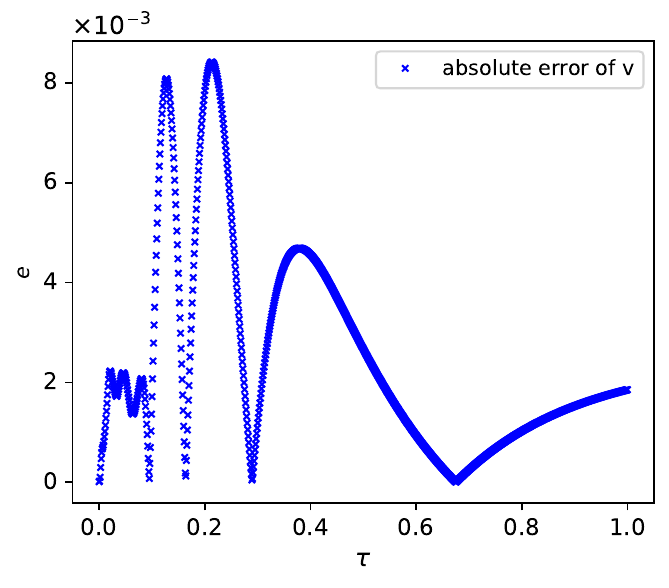}
        }
        \subfigure[absolute error of reference and vanilla PINN solutions of $z$]
    {\includegraphics[width=0.28\textwidth]{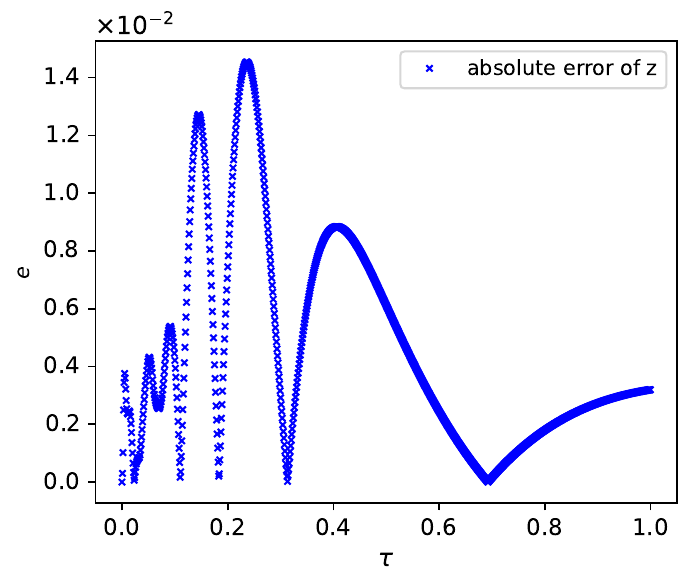}}
        \subfigure[absolute error of reference and vanilla PINN solutions of $w$]
        {
        \includegraphics[width=0.28\textwidth]{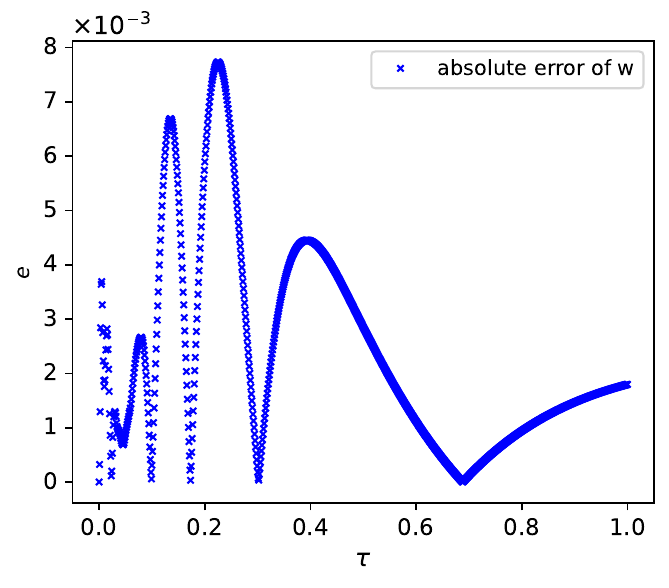}
        }
        \caption{Results for Example \ref{exm:FNsys}  when $\epsilon_1=10^{-1}, \epsilon_2=\epsilon_3=10^{-2}$.}
 \label{fig:results_FHN_1e-2}
        
\end{figure}

\begin{figure}[!htb]
	\centering
	\subfigure[reference and 2SNN solutions \quad of $v$]
      {\includegraphics[width=0.30\textwidth]{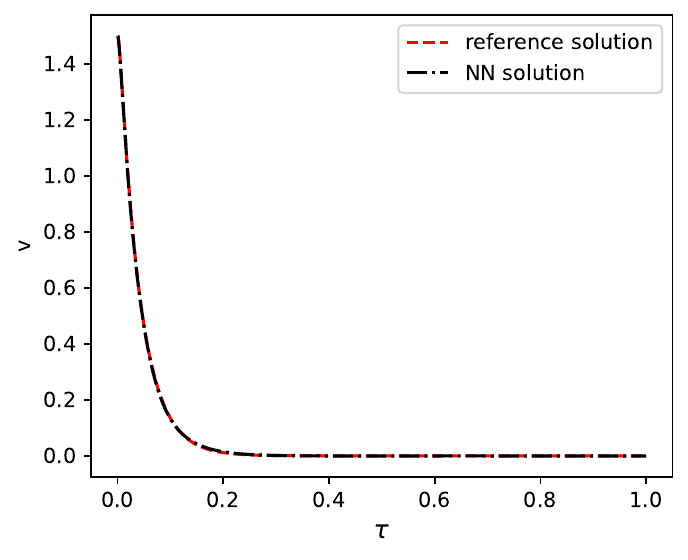}
        }
        \subfigure[reference and 2SNN solutions \quad of $z$]
        {\includegraphics[width=0.30\textwidth]{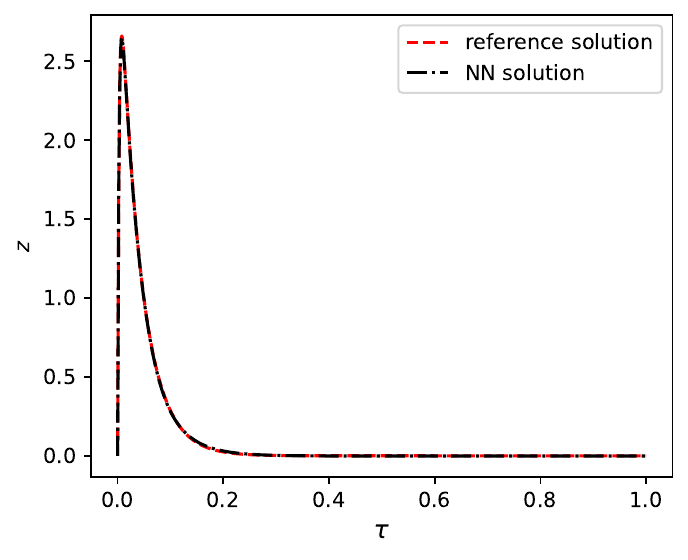}}
        \subfigure[reference and 2SNN solutions \quad of $w$]
        {\includegraphics[width=0.30\textwidth]{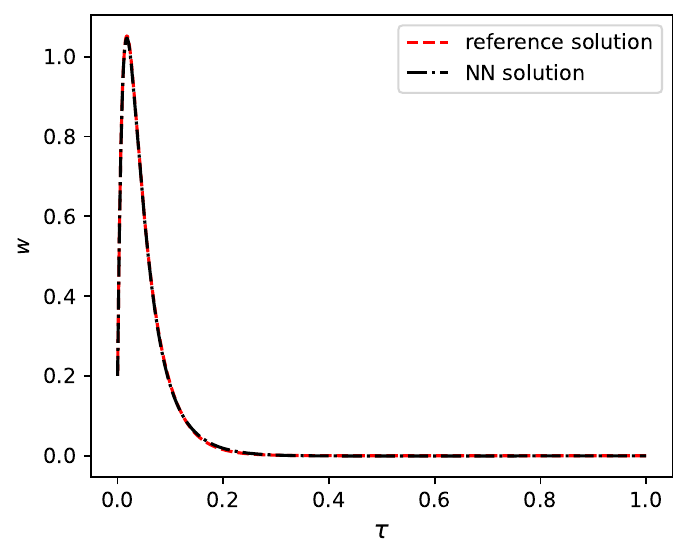}
        }

        \subfigure[absolute error of reference and 2SNN solutions of $v$]
      {\includegraphics[width=0.30\textwidth]{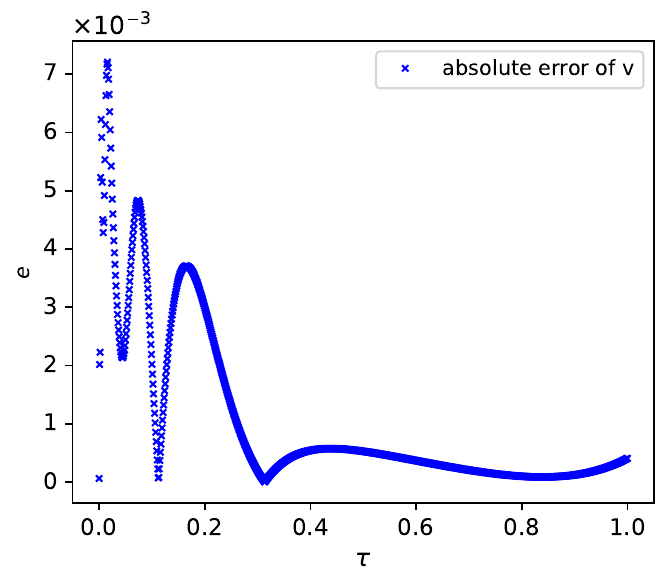}
       
        }
        \subfigure[absolute error of reference and 2SNN solutions of $z$]
        {\includegraphics[width=0.30\textwidth]{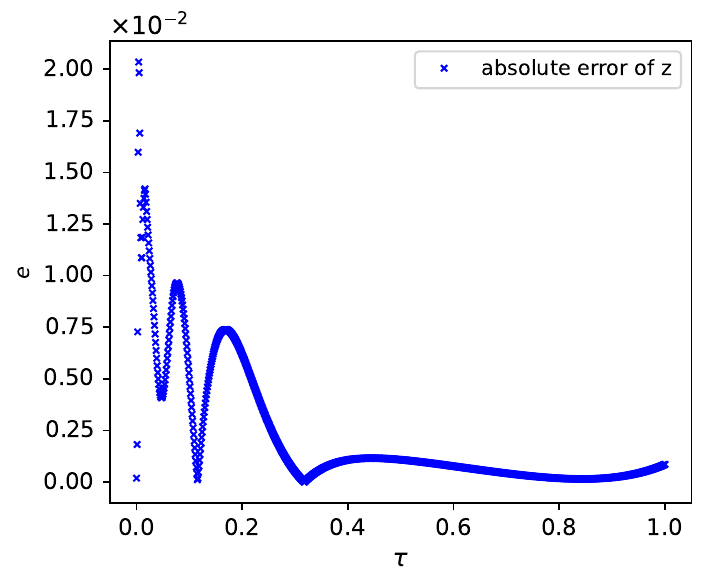}}
        \subfigure[absolute error of reference and 2SNN solutions of $w$]
        {
        \includegraphics[width=0.30\textwidth]{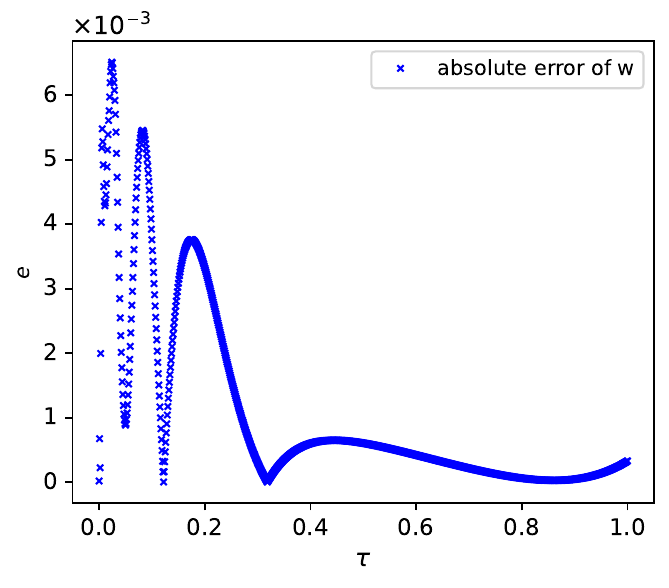}
        }
        \subfigure[absolute error of reference and vanilla PINN solutions of $v$]
      {\includegraphics[width=0.30\textwidth]{figs/exm_FHN/FNsys_abserr_v_1s1e-2.pdf}
        }
        \subfigure[absolute error of reference and vanilla PINN solutions of $z$]
    {\includegraphics[width=0.31\textwidth]{figs/exm_FHN/FNsys_abserr_z_1s1e-2.pdf}}
        \subfigure[absolute error of reference and vanilla PINN solutions of $w$]
        {
        \includegraphics[width=0.30\textwidth]{figs/exm_FHN/FNsys_abserr_w_1s1e-2.pdf}
        }
        \caption{Results for Example \ref{exm:FNsys}  when $\epsilon_1=10^{-1}, \epsilon_2=10^{-2}/8, \epsilon_3=10^{-2}$ \emph{with} successive training (trained from the same values of $\epsilon_1$ and $\epsilon_3$, but $\epsilon_2=10^{-2}/4$). }
 \label{fig:results_FHN_1e-2_8}
\end{figure}

\begin{figure}[!htb]
	\centering
	\subfigure[reference and 2SNN solutions \quad of $v$]
      {\includegraphics[width=0.30\textwidth]{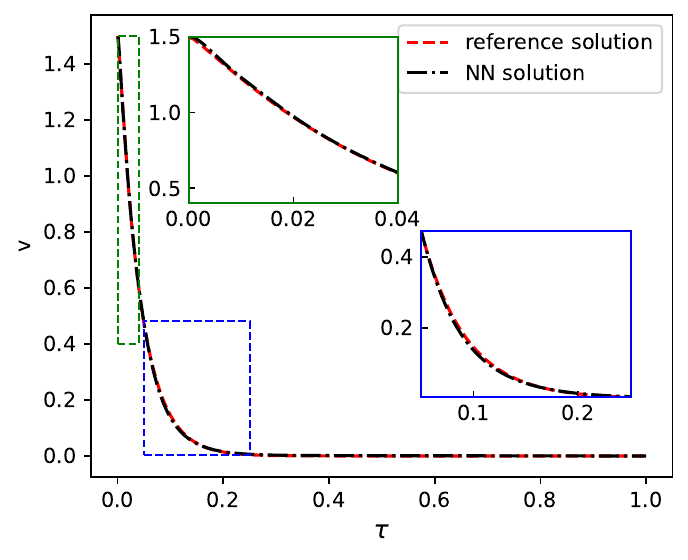}
        }
        \subfigure[reference and 2SNN solutions \quad of $z$]
        {\includegraphics[width=0.30\textwidth]{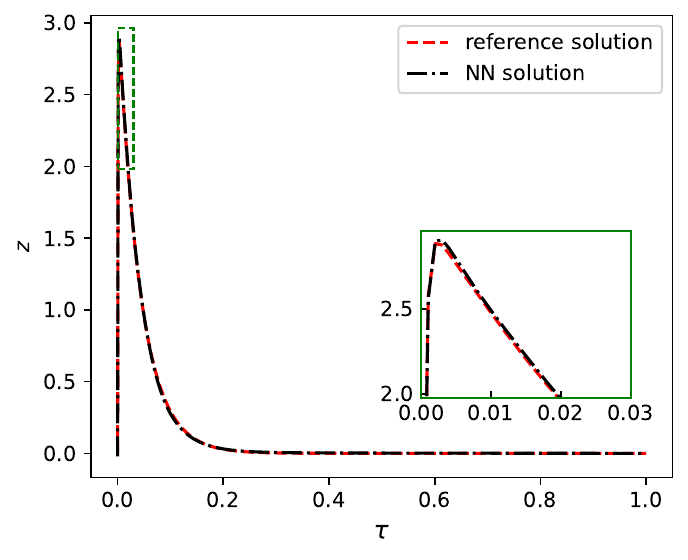}}
        \subfigure[reference and 2SNN solutions \quad of $w$]
        {\includegraphics[width=0.30\textwidth]{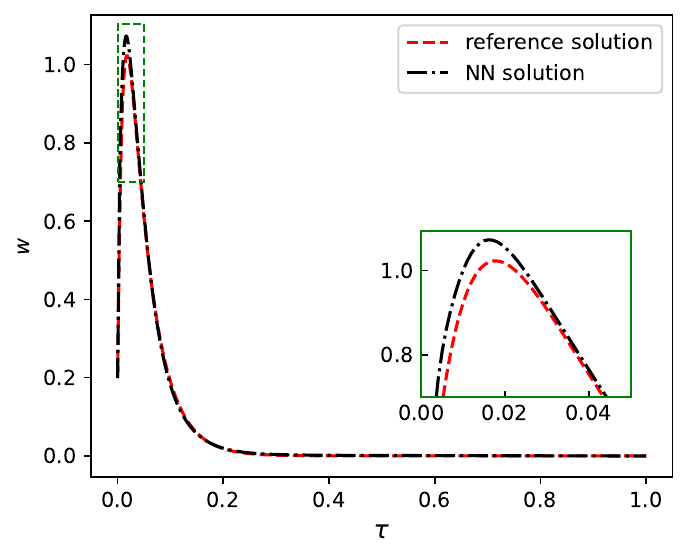}
        }

        \subfigure[reference and vanilla PINN solutions of $v$]
      {\includegraphics[width=0.30\textwidth]{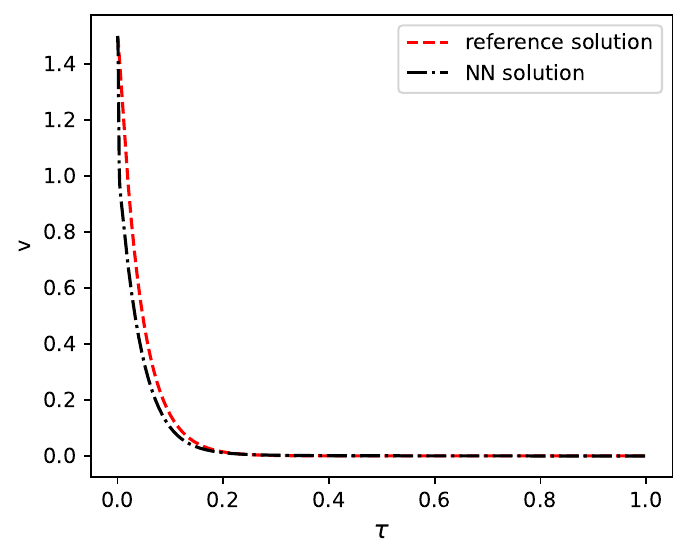}
       
        }
        \subfigure[reference and vanilla PINN solutions of $z$]
        {\includegraphics[width=0.30\textwidth]{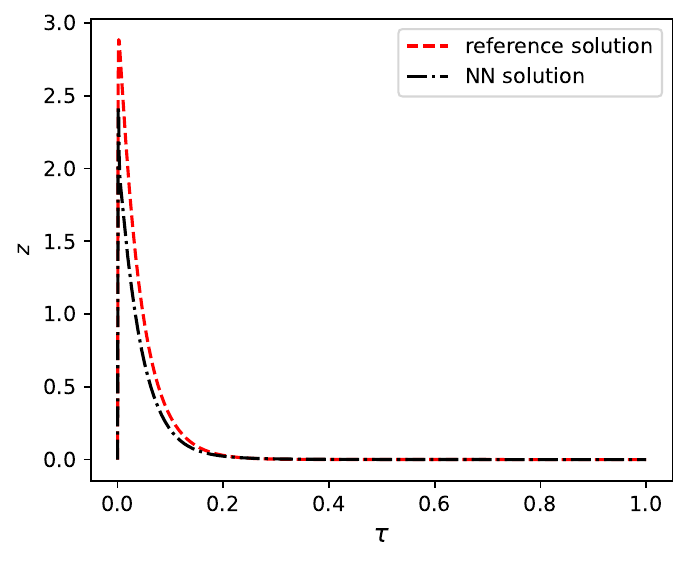}}
        \subfigure[reference and vanilla PINN solutions of $w$]
        {
        \includegraphics[width=0.30\textwidth]{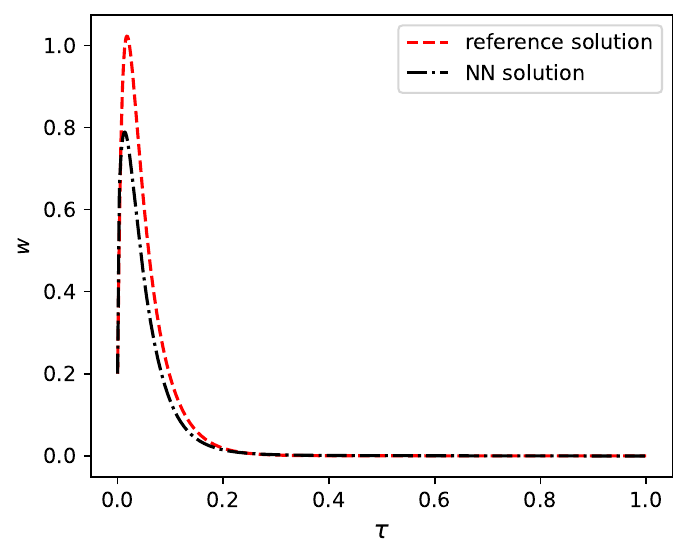}
        }
        \caption{Results for Example~\ref{exm:FNsys} obtained \emph{with} successive training pertinent to Table \ref{tab:FNsys_para}, with $\epsilon_1=10^{-1}$, $\epsilon_2=2.5\times 10^{-4}$, and $\epsilon_3=10^{-2}$.}

 \label{fig:results_FHN_1e-3_4_suc}
        
\end{figure}

\begin{figure}[!htb]
	\centering
     \subfigure[trained from $\epsilon_2=10^{-2}/8$ to $10^{-2}/16$.]{
        \includegraphics[width=0.28\textwidth]{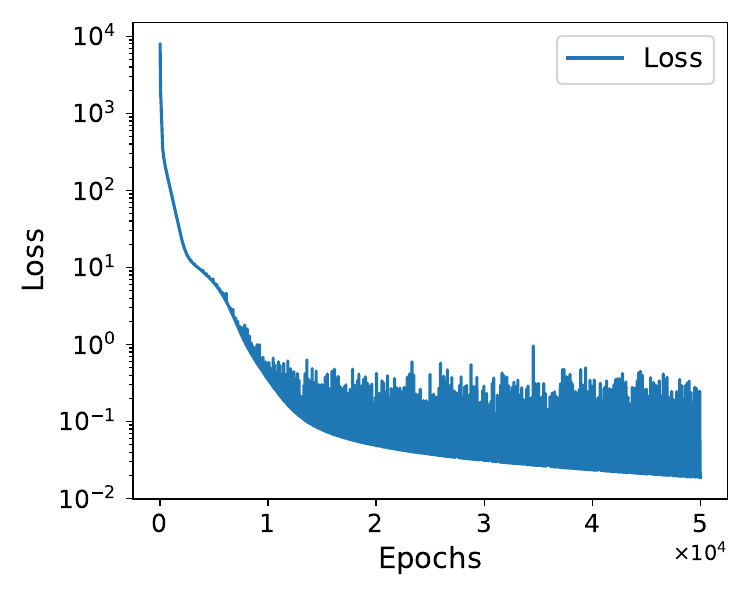}
        }
        \subfigure[trained from $\epsilon_2=10^{-2}/16$ to $10^{-2}/32$]{
        \includegraphics[width=0.28\textwidth]{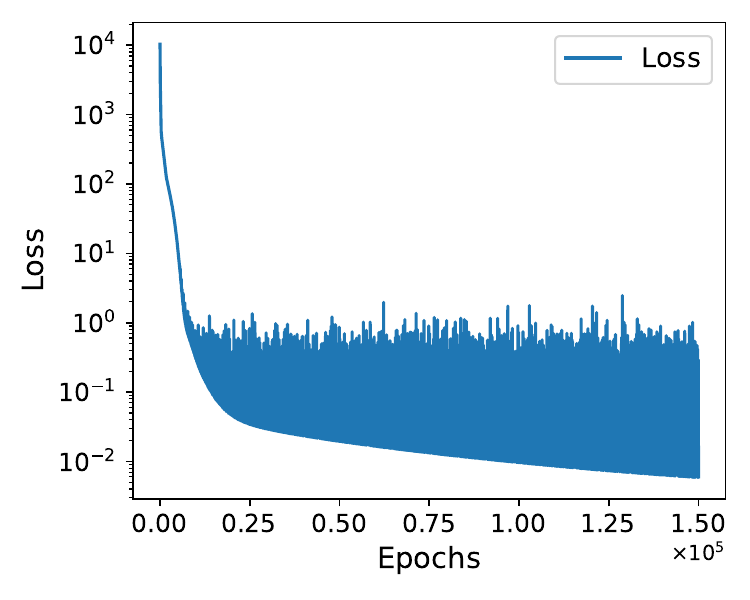}
        }
        \subfigure[trained from $\epsilon_2=10^{-2}/32$ to $2.5\times 10^{-4}$]{
        \includegraphics[width=0.28\textwidth]{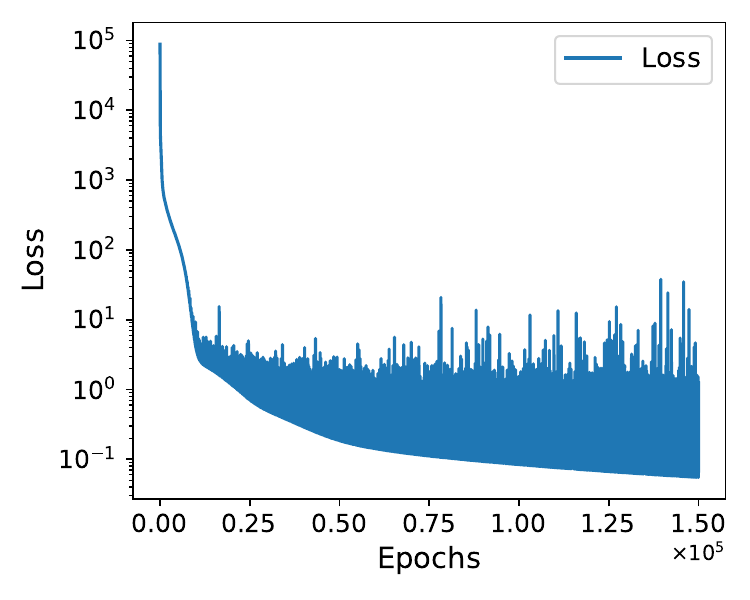}
        }
        \caption{ Loss history for Examples \ref{exm:FNsys} using successive training strategy with parameters in Table \ref{tab:FNsys_para}.
}
		 \label{fig:exmFHN_MSNN_loss}
\end{figure}

\begin{table}[!htb]
\centering
 \begin{adjustbox}{width=0.86\textwidth}
\begin{tabular}{c c| c c c | c c}
\hline
$\epsilon_2$ & $10^{-2}$ & $10^{-2}/4$ & $10^{-2}/8$ & $10^{-2}/16$ & $10^{-2}/32$ & $2.5\times 10^{-4}$ \\
$\epsilon$ & $2.15\times10^{-2}$ & $1.36\times10^{-2}$ & $1.08\times10^{-2}$ & $8.55\times 10^{-3}$ & $6.79\times 10^{-3}$ & $6.30\times 10^{-3}$ \\
\hline
$\alpha$ & 1000 &\multicolumn{3}{c|}{$1000$} & \multicolumn{2}{c}{$1000$} \\
$N_c$ & 450 & \multicolumn{3}{c|}{$450$}& \multicolumn{2}{c}{450} \\
LR & P-S & \multicolumn{3}{c|}{$10^{-4}$} &\multicolumn{2}{c}{$10^{-4}$} \\
iterations & $5\times 10^4$ & \multicolumn{3}{c|}{$5\times10^4$} & \multicolumn{2}{c}{$1.5\times10^5$} \\
\hline
\end{tabular}
\end{adjustbox}
\vskip 3pt
\caption{Parameters in the loss function \eqref{eq:loss-general} and hyper-parameters of the successive training for Example~\ref{exm:FNsys}. LR stands for learning rate, P-S is the piecewise constant scheduler in Table~\ref{tab:lr_pt}.}
\label{tab:FNsys_para}
\end{table}
We report the results of successive training in Figures \ref{fig:results_FHN_1e-2}, \ref{fig:results_FHN_1e-2_8}, and \ref{fig:results_FHN_1e-3_4_suc} and Table \ref{tab:FHN_suctrain_intermediate}, with the related training loss history in Figure \ref{fig:exmFHN_MSNN_loss}. In the initial step of the successive training ($\epsilon_2=10^{-2}$), the 2SNN solution agrees with the reference solution well, as indicated in Figure \ref{fig:results_FHN_1e-2}. For the intermediate stages when $\epsilon_2=10^{-2}/2^i$ $(i=2,3,4,5)$, we compute the errors between 2SNN and reference solutions in $l^2$ and $l^\infty$ norms and record them in Table \ref{tab:FHN_suctrain_intermediate}. Under those setups, the neural network and reference solution comparisons resemble Figure \ref{fig:results_FHN_1e-2_8}, indicating a reasonably accurate match.  
Those figures and Table \ref{tab:FHN_suctrain_intermediate} indicate that during the initial and intermediate stages of successive training, the 2SNN and vanilla PINN provide comparable accuracy.

However, when $\epsilon_2$ decreases to $10^{-3}/4$ (a high-contrast regime with $\epsilon_1/\epsilon_2 = 400$), the 2SNN \emph{significantly} outperforms the vanilla PINN, as shown in Figure \ref{fig:results_FHN_1e-3_4_suc}. Figure \ref{fig:results_FHN_1e-3_4_suc} (a)-(c) show that the 2SNN solutions agree reasonably well with the reference solutions, whereas the corresponding Figure \ref{fig:results_FHN_1e-3_4_suc} (d)-(f) for the vanilla PINN exhibit substantial deviations from the reference ones, especially near the initial layers of the solutions for $z$ and $w$.

\begin{figure}[!htb]
	\centering
	\subfigure[reference and 2SNN solutions \quad of $v$]
      {\includegraphics[width=0.30\textwidth]{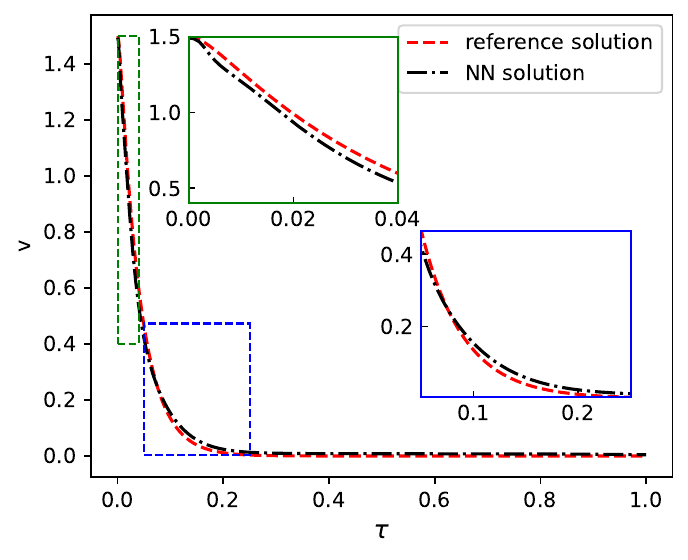}
        }
        \subfigure[reference and 2SNN solutions \quad of $z$]
        {\includegraphics[width=0.30\textwidth]{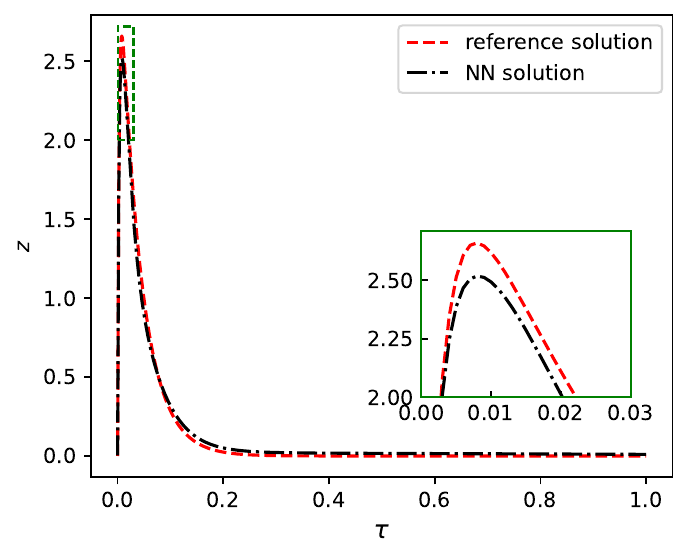}}
        \subfigure[reference and 2SNN solutions \quad of $w$]
        {\includegraphics[width=0.30\textwidth]{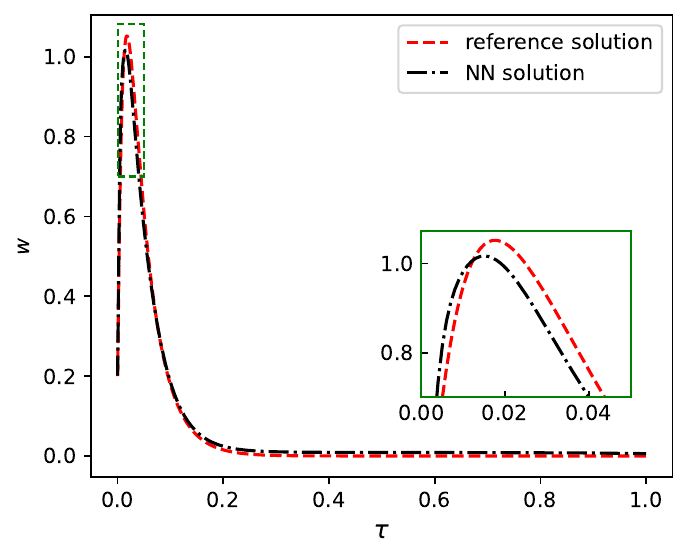}
        }

        \subfigure[reference and vanilla PINN solutions of $v$]
      {\includegraphics[width=0.30\textwidth]{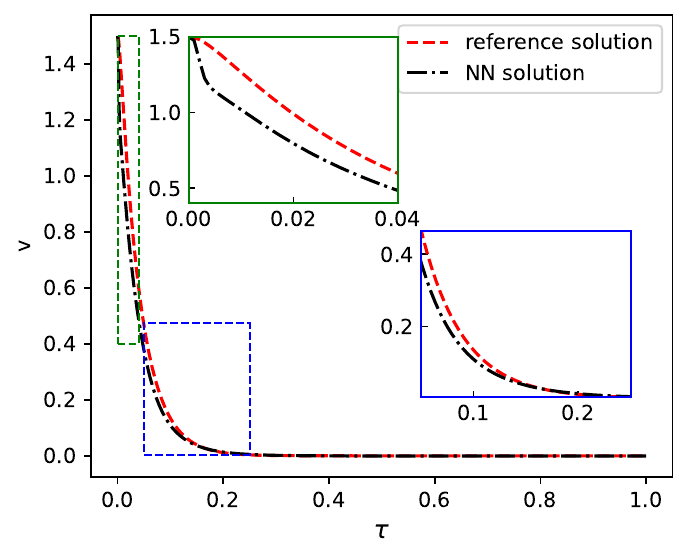}
       
        }
        \subfigure[reference and vanilla PINN solutions of $z$]
        {\includegraphics[width=0.30\textwidth]{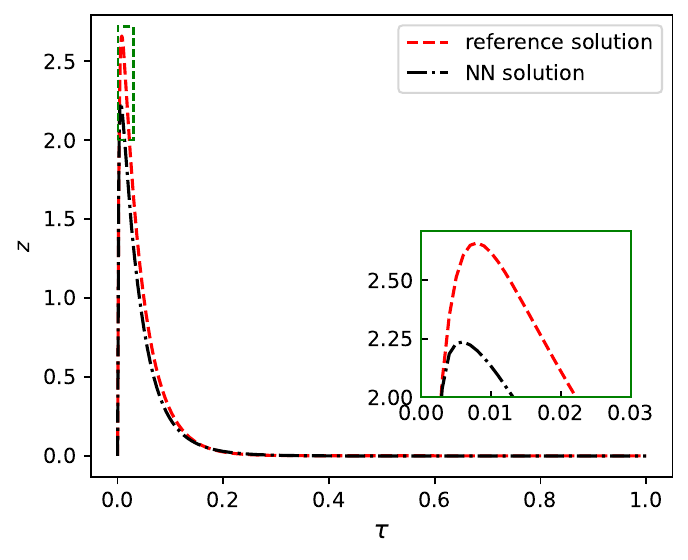}}
        \subfigure[reference and vanilla PINN solutions of $w$]
        {
        \includegraphics[width=0.30\textwidth]{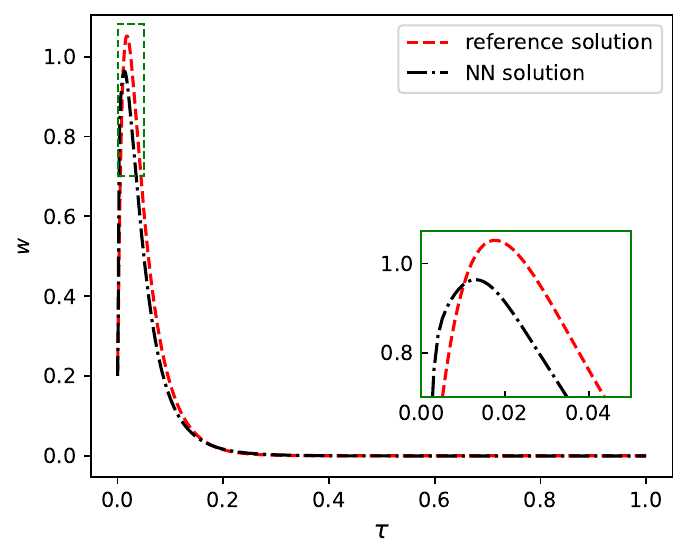}
        }
        \caption{Results for Example~\ref{exm:FNsys} obtained \emph{without} successive training, with $\epsilon_1=10^{-1}$, $\epsilon_2=10^{-2}/8$, and $\epsilon_3=10^{-2}$.}

 \label{fig:results_FHN_1e-2_8_direct}
        
\end{figure}

We also report the results of directly solving for the case $\epsilon_1=10^{-1}, \epsilon_3=10^{-2}, \epsilon_2=10^{-2}/8$ in Figure \ref{fig:results_FHN_1e-2_8_direct} and Table \ref{tab:error_FHN_1e-2_8_direct} \emph{without} the successive training. As indicated in Figure \ref{fig:results_FHN_1e-2_8_direct}, with the zoom-in views using green boxes for solutions around initial layers, the accuracy of 2SNN solutions (shown in Figure \ref{fig:results_FHN_1e-2_8_direct} (a)-(c)) outperforms the vanilla PINN ones (Figure \ref{fig:results_FHN_1e-2_8_direct} (d)-(f)) around the initial layers. This is especially true for the components $z$ and $w$, where vanilla PINN solutions deviate significantly from the reference, whereas 2SNN solutions provide reasonable predictions. Meanwhile, we observe that the vanilla PINN solution slightly outperforms the 2SNN solutions for tails away from the initial layer, especially from the zoom-in views in blue boxes in Figure \ref{fig:results_FHN_1e-2_8_direct} (a) and (d). We also record the errors for $v$, $z$, and $w$ in Table \ref{tab:error_FHN_1e-2_8_direct}. The table shows that the 2SNN solutions outperform vanilla PINN solutions across these error metrics. 

A comparison in errors from Tables \ref{tab:error_FHN_1e-2_8_direct} and \ref{tab:FHN_suctrain_intermediate}  shows that when $\epsilon_2 = 10^{-2}/8$ (as highlighted in the tables), curriculum learning (successive training) leads to improved accuracy. This improvement suggests enhanced training stability through better parameter initialization at each stage of Step 2 in Algorithm \ref{alg:succesive-training}.

To summarize, Example \ref{exm:Robertson} highlights the importance of a suitable choice of the effective parameter, as stated in Claim C, while Example \ref{exm:FNsys} demonstrates that curriculum learning improves stability, as stated in Claim D. Both examples confirm accurate predictions under regimes of high-contrast small parameters, in line with Claim B. Also, both examples demonstrate Claim A, showing that 2SNN properly resolves the effects of coupling. 

\begin{table}[!htp]
    \centering
    \vspace{2mm}
    \renewcommand{\arraystretch}{1.2}
    \begin{adjustbox}{width=0.89\textwidth}
    \begin{tabular}{c|cccccc}
        \toprule
        \textbf{Method}
        & $\lVert e_v \rVert_{l_\infty}$
        & $\lVert e_v \rVert_{l_2}$
        & $\lVert e_z \rVert_{l_\infty}$
        & $\lVert e_z \rVert_{l_2}$
        & $\lVert e_w \rVert_{l_\infty}$
        & $\lVert e_w \rVert_{l_2}$ \\
        \midrule
        2SNN
        & $\mathbf{7.38\times 10^{-2}}$ & $\mathbf{1.67\times 10^{-2}}$ & $\mathbf{1.44\times 10^{-1}}$ & $\mathbf{3.32\times 10^{-2}}$ & $\mathbf{1.16\times 10^{-1}}$ & $\mathbf{1.62\times 10^{-2}}$ \\
        vanilla PINN
& $2.74\times 10^{-1}$ & $4.19\times 10^{-2}$ & $4.87\times 10^{-1}$ & $8.02\times 10^{-2}$ & $2.26\times 10^{-1}$ & $3.35\times 10^{-2}$ \\        \bottomrule
    \end{tabular}
    \end{adjustbox}
    \vskip 3pt
    \caption{Errors of 2SNN and vanilla PINN solutions Example~\ref{exm:FNsys} when $\epsilon_1=10^{-1}, \epsilon_2=10^{-2}/8, \epsilon_3=10^{-2}$, \emph{without} successive training.}
    \label{tab:error_FHN_1e-2_8_direct}
\end{table}

\begin{table}[!htp]
\centering
\renewcommand{\arraystretch}{1.2}
\begin{adjustbox}{width=0.89\textwidth}
\begin{tabular}{c c c c c c c}
\toprule
\textbf{Method} / $\epsilon_2$
& $\|e_v\|_{l_\infty}$
& $\|e_v\|_{l_2}$
& $\|e_z\|_{l_\infty}$
& $\|e_z\|_{l_2}$
& $\|e_w\|_{l_\infty}$
& $\|e_w\|_{l_2}$ \\
\midrule
\multicolumn{7}{l}{\textbf{2SNN}} \\
$\epsilon_2 = 10^{-2}/4$    & $5.85\times 10^{-3}$ & $1.97\times 10^{-3}$ & ${1.18\times 10^{-2}}$ & $3.96\times 10^{-3}$ & $5.65\times 10^{-3}$ & $1.89\times 10^{-3}$
 \\
$\mathbf{\epsilon_2 = 10^{-2}/8}$    & $\mathbf{7.21\times 10^{-3}}$ & $\mathbf{1.73\times 10^{-3}}$ & $\mathbf{{2.03\times 10^{-2}}}$ & $\mathbf{3.67\times 10^{-3}}$ &$\mathbf{6.52\times 10^{-3}}$ & $\mathbf{1.80\times 10^{-3}}$ \\
$\epsilon_2 = 10^{-2}/16$  &$1.95\times 10^{-2}$ & $3.14\times 10^{-3}$ & ${3.56\times 10^{-2}}$ & $6.17\times 10^{-3}$ & $1.33\times 10^{-2}$ & $2.85\times 10^{-3}$
 \\
$\epsilon_2 = 10^{-2}/32$   & $3.13\times 10^{-2}$ & $4.69\times 10^{-3}$ & ${5.79\times 10^{-2}}$ & $9.14\times 10^{-3}$ & $3.05\times 10^{-2}$ & $5.55\times 10^{-3}$

 \\
\midrule
\multicolumn{7}{l}{\textbf{vanilla PINN}} \\

$\epsilon_2 = 10^{-2}/4$    & $2.06\times 10^{-3}$ & $5.43\times 10^{-4}$ & $4.90\times 10^{-3}$ & $1.11\times 10^{-3}$ & $5.40\times 10^{-3}$ & $8.21\times 10^{-4}$
 \\
$\epsilon_2 = 10^{-2}/8$    & $1.06\times 10^{-2}$ & $2.10\times 10^{-3}$ & $2.76\times 10^{-2}$ & $4.52\times 10^{-3}$ & $9.53\times 10^{-3}$ & $2.11\times 10^{-3}$
 \\
$\epsilon_2 = 10^{-2}/16$   & $6.30\times 10^{-3}$ & $1.86\times 10^{-3}$ & $4.42\times 10^{-2}$ & $4.32\times 10^{-3}$ & $6.79\times 10^{-3}$ & $1.86\times 10^{-3}$
 \\
$\epsilon_2 = 10^{-2}/32$   & $8.94\times 10^{-3}$ & $1.81\times 10^{-3}$ & $6.06\times 10^{-2}$ & $4.44\times 10^{-3}$ & $5.08\times 10^{-3}$ & $1.73\times 10^{-3}$
 \\
\bottomrule
\end{tabular}
\end{adjustbox}
\vskip 3pt
\caption{Errors of 2SNN and vanilla PINN solutions for Example \ref{exm:FNsys} under different choices of $\epsilon_2$, with $\epsilon_1=10^{-1}$ and $\epsilon_3=10^{-2}$, \emph{with} successive training.}
\label{tab:FHN_suctrain_intermediate}
\end{table}

\section{Conclusion and Discussion}\label{sec:conclusion}
In this study, we apply the  two-scale neural network framework to dynamical systems with multiple small parameters using a new scale parameter from the geometric mean of all small parameters. We have the following key findings in the aspects of methodology, theoretical justification, and empirical validation:

(i) \emph{Methodology}: The network input is augmented with a scale-aware feature based on a locally stretched coordinate and a single effective scale parameter representing the multiple small parameters in the model, thereby intrinsically accommodating multiple scales in multi-small-parameter systems in a streamlined manner. A curriculum learning scheme is applied to stabilize training.
(ii) \emph{Theoretical justification}: We provide theoretical heuristics for the proposed network design for singularly perturbed systems with multiple small parameters, and explain the rationale for selecting the geometric mean as an effective scale parameter to represent these small parameters.
(iii) \emph{Empirical validation}: Numerical experiments on a range of dynamical systems demonstrate the effectiveness and accuracy of the proposed framework in capturing sharp solution transitions induced by small parameters. 
%
%

Limitations persist in regimes with extreme stiffness or scale hierarchies that are poorly separated or unknown, where additional specialized treatments may be necessary to ensure robustness and applicability. When extreme stiffness arises from high-contrast model parameters, using the geometric mean as the effective (aggregated) small parameter may under-resolve the underlying multiscale structure. Such challenges may be alleviated by using high-order optimizers, e.g., in \cite{an2026lightweight,jnini2026curvature,Urban2025unveiling}. 
%
For problems with high-contrast parameters, the proposed approach can serve as a reasonable initialization and may be integrated into more sophisticated network frameworks (e.g., multi-stage neural networks \cite{LaiMultistage2024}) or training patterns to achieve the desired accuracy, which will be explored in future work.

We expect that the present two-scale neural network framework can facilitate learning operators arising from singularly perturbed systems with multiple parameters, which generalizes
operator learning in \cite{PASNet2025} for scalar singularly perturbed partial differential equations. The proposed neural network framework can also be applied in mechanistic models with data to investigate large-scale fast-slow dynamical systems, such as those in the health sciences \cite{biosecureonehealth_2025, Abebiophysismodel2018}.

\section* {Acknowledgment}
We thank Mr. Truong Hoang Nhan Pham, formerly of the University of Missouri-Kansas City, for assistance in preparing a portion of the preliminary materials related to this work.

\section*{Funding statement}
There is no funding for this research.

\section*{Declaration of competing interest}
The authors declare that they have no known competing financial interests or personal relationships that could have influenced the work reported in this manuscript. This research was conducted independently, without any commercial or financial support that could be construed as a potential conflict of interest.
All authors have reviewed and approved the final version of the manuscript and agree with its submission. 
The authors affirm that the work represents their original research, has not been published previously, and is not under consideration elsewhere.
\section*{Data availability statements}
The code and data that support the findings of this study are available on request from the corresponding author.
\section*{Declaration of generative AI and AI-assisted technologies in the manuscript preparation}
During the preparation of this manuscript, the authors used generative AI-based tools solely for language polishing, proofreading, and improving the clarity and readability of the text.
The authors affirm that they take full responsibility for the content of the manuscript. All scientific ideas, methodologies, algorithms, experimental designs, benchmark tests, and conclusions presented in this work were conceived, developed, and validated by the authors. The use of AI tools did not influence the scientific content or the interpretation of the results.

\bibliographystyle{plain}
\bibliography{singptb-ref,nn,MsNN}

\end{document}